\documentclass[a4paper,11pt]{amsart}

	\usepackage[utf8]{inputenc}
	\usepackage[T1]{fontenc}
	\usepackage[english]{babel}
	\usepackage[protrusion=true,expansion=true,kerning, spacing]{microtype}
	\microtypecontext{spacing=english}
	\usepackage{CJKutf8} 
	\usepackage{cmap} 

	\usepackage{xspace}
	
	\usepackage[textwidth=16cm,hcentering,heightrounded,foot=1cm]{geometry}
	\usepackage{caption}
 \usepackage{titlesec}
	\titleformat{name=\section}[hang]{\scshape \large \centering}{\thesection.~}{0em}{}[]

	\titleformat{\subsection}[hang]{\bfseries}{\thesubsection.~}{0em}{}[]

	\titleformat{\subsubsection}[runin]{\bfseries}{\thesubsubsection.}{.2em}{}[]

 \usepackage{titletoc}
\titlecontents{section}[2pc]{\addvspace{1pc}}
{\contentslabel[\textsc{\thecontentslabel}]{1pc}}
{}{\dotfill\contentspage}
[\addvspace{2pt}]

\titlecontents{subsection}[4pc]{}
{\contentslabel[\subsectionname\thecontentslabel]{2pc}}
{}{\hfill\contentspage}
[]

\usepackage{csquotes}
\usepackage[
bibencoding=utf8,%
giveninits=true,
bibstyle=ieee-alphabetic,
citestyle=alphabetic,%
sorting=nyt,%
language=french,%
isbn=false,%
doi=false,%
eprint=true,
backref=true, 
hyperref=true,
url=true,
minbibnames=4,
maxbibnames=4,
dashed=true,
backend=biber,
]{biblatex}
	\usepackage{amsthm,amssymb,amsmath}
	\usepackage{amssymb,amsfonts}
	\usepackage[fixamsmath]{mathtools} 
	\usepackage{accents} 
	\usepackage{tikz-cd}
		\numberwithin{equation}{subsection}

	\usepackage{graphics,graphicx}
	\usepackage{color}
	\usepackage{tikz}
	\usepackage{caption} 
	\usepackage{subcaption} 
	\usepackage{wrapfig}

\theoremstyle{definition}
	\newtheorem{defin}{Definition}[section]
	\newtheorem{ex}[defin]{Example}

\theoremstyle{definition}

\theoremstyle{plain}
	\newtheorem{theo}[defin]{Theorem}
	\newtheorem{prop}[defin]{Proposition}
	\newtheorem{cor}[defin]{Corollary}
	\newtheorem{lem}[defin]{Lemma}

\newcommand{\C}{\mathbf{C}}

\newcommand{\N}{\mathbf{N}}
\newcommand{\R}{\mathbf{R}}
\newcommand{\Q}{\mathbf{Q}}
\newcommand{\Z}{\mathbf{Z}}

\newcommand{\F}{\mathbf{F}}

\newcommand{\G}{\mathbf{G}}

\newcommand{\HH}{\mathbf{H}}

\newcommand{\OO} {{\mathcal O}} 

\newcommand{\End} {\operatorname{End} }

\newcommand{\Card} {\operatorname{Card} }
\newcommand{\Aut} {\operatorname{Aut} }

\newcommand{\GL}{\operatorname{GL}}

\newcommand{\Sp}{\operatorname{Sp}}

\usepackage[pdftex,
hidelinks,
hyperindex=true, bookmarks=true, bookmarksopen=true, bookmarksnumbered=true,
pdfauthor={Séverin Philip},
pdftitle={The (p,t,a)-inertial groups as finite monodromy groups},
pdfsubject={},
pdfstartview=FitH
]{hyperref}

\begin{document}
\author{Séverin Philip}
\address{Department of Mathematics, Stockholms universitet
SE-106 91 Stockholm, Sweden}
	\email{severin.philip@math.su.se}
\thanks{The author is grateful to J.~Bergström, G.~Rémond and A.~Tamagawa for the fruitful discussions about this project.}
\title{The $(p,t,a)$-inertial groups as finite monodromy groups \\[0.5em]
\textnormal{ Version of \today}}

\begin{abstract}
    Silverberg and Zarhin introduced the notion of a $(p,t,a)$-inertial group in the hope of having a group theoretic characterization of the finite groups that appear as finite monodromy groups -- the groups that represent the local obstruction to semi-stable reduction -- of abelian varieties in fixed dimension $t+a$. In this text, we provide a positive answer to their question, that is, every $(p,t,a)$-inertial group is the finite monodromy group of an abelian variety in dimension $t+a$. To prove this, we show a structure theorem on the rational group algebra $\Q[G]$ of ramification groups, refining a theorem of Serre and generalizing results on $p$-groups of Roquette and Ford.
\end{abstract}

\maketitle

\tableofcontents

\section{Introduction}

\stepcounter{subsection}

\subsubsection{} The finite monodromy groups of abelian varieties have been introduced by Grothendieck in \cite{sga} exposé IX. They represent the local obstruction to semi-stable reduction. Silverberg and Zarhin studied these groups in \cite{SZ98, SZ05}, and the author, in order to give an effective version of Grothendieck's semi-stable reduction theorem, in \cite{Ph221, Phi222, Ph4}.

\medskip

For a fixed natural integer $g\geq 1$, the list of finite groups which can be realized as finite monodromy groups of some abelian variety of dimension $g$ is not known. An attempt to provide this list is made by the notion of $(p,t,a)$-inertial groups introduced in \cite{SZ05}. 

\begin{defin}[\cite{SZ05}] \label{def:pta} Let $p$ be a prime number or $p=0$ and $t$, $a$ positive integers. A finite group $G$ is said to be $(p,t,a)$-inertial if it satisfies the two following conditions :
\begin{itemize}
    \item[$(i)$] If $p=0$ then $G$ is cyclic otherwise $G$ is a semi-product $\Gamma_p\rtimes \Z/n\Z$ with $\Gamma_p$ a $p$-group and $n$ an integer prime to $p$.
    \item[$(ii)$] For all primes $\ell\neq p$ there is an injection
    \[
    \iota_{\ell}\colon G\hookrightarrow \GL_t(\Z) \times \Sp_{2a} (\Q_{\ell})
    \]
    
\end{itemize}
    such that the projection map onto the first factor is independent of $\ell$ and the characteristic polynomial of the projection of any element onto the second factor has integer coefficients independent of $\ell$. 
\end{defin}

It follows from the definition that the character $\chi_{\ell}$ of the $\Q_{\ell}$ representation given by the projection on the second factor has integer values and is independent of $\ell$. Note that from the main theorem of \cite{Mau68}, condition $(i)$ is equivalent to the fact that $G$ is a quotient of the absolute inertia group of some local fields of equal and mixed characteristic -- a ramification group at $p$. 

\medskip

It is known -- see part $(i)$ of Theorem~5.2 of \cite{SZ98} -- that finite monodromy groups are $(p,t,a)$-inertial. The question is now whether the converse is true, that is, does Definition~\ref{def:pta} characterize the finite groups that appear as finite monodromy groups in dimension $t+a$ over number fields with conditions on the toric and abelian ranks of the reduction with respect to $t$ and $a$. A precise statement is given by question 1.13 of \textit{loc. cit.} Specifically, the question is to realize a $(p,t,a)$-inertial group $G$ as the finite monodromy group of an abelian variety $A$ of dimension $g=t+a$ at a place of residue characteristic $p$. In the same paper they show that the set of $(p,t,a)$-inertial groups with $t+a=2$ are realized as finite monodromy groups of abelian surfaces over local fields of equal characteristic $p$. Further, Chrétien and Matignon in \cite{CM13} have shown that this list is also realized with abelian surfaces over number fields by an ad hoc construction for the last unrealized $(2,0,2)$-inertial group. 

\subsubsection{} In this text we give a positive answer to their question.

\begin{theo}[Corollary~\ref{cor:realmonod}]
     Let $G$ be a finite group. Then $G$ is $(p,t,a)$-inertial if and only if $G$ is the finite monodromy group of an abelian variety $A$ of dimension $t+a$ over a number field $K$ at a place $v$ of residue characteristic $p$ and $A$ is such that, for all extension $L/K$ for which $A_L$ has semi-stable reduction at the places of $L$ above $v$, the reduction of $A_L$ at the places above $v$ has toric rank $t$ and abelian rank $a$. 
\end{theo}

We actually show more in Theorem~\ref{theo:realmonod}. Indeed, we show that for $G$ a ramification group and every family of maps as in $(ii)$ of Definition~\ref{def:pta} there is a realization of $G$ as a finite monodromy group of an abelian variety of dimension $t+a$ over a $p$-adic field such that we recover that family of maps, up to isomorphism, by taking the $\ell$-adic representation of the semi-stable reduction of $A$. In order to prove this result we deal with the rational group algebra $\Q[G]$ of such groups and give precision to a structure theorem of Serre, Théorème~3 of \cite{Se60}. 

\medskip

We end the last section by a series of example, showcasing how the results in this paper are used to compute the values for which a ramification group is $(p,t,a)$-inertial in a simple and efficient way.

\begin{theo}[Theorem~\ref{theo:structureQG}] \label{theo:centresalgebresintro}
   Let $G$ be a ramification group at $p$. The rational group algebra $\Q[G]$ has the following properties. 
    \begin{itemize}
        \item[$(i)$] The algebra $\Q[G]$ is quasi-split outside $p$.
        \item[$(ii)$] If $E\subset \Q[G]$ is a simple factor of $\Q[G]$ then the center $Z(E)$ of $E$ is a CM field or a subfield of $\Q(\mu_{p^{\infty}})$.
    \end{itemize}
    Furthermore, if $\Card G$ is odd then no simple factor $E$ of $\Q[G]$ is such that $E\otimes_{\Q} \R$ is a matrix algebra over the quaternions.  
\end{theo}

Part $(i)$ is the statement of Théorème~3 of \cite{Se60}. We deal with the possible centers of the simple factors of $\Q[G]$. The main novelty obtained from $(ii)$ is that in the event of a simple factor $E$ of $\Q[G]$ being a type III algebra in Albert's classification, then there is a totally real subfield $F$ of $\Q(\mu_p^{\infty})$ such that we have $E\simeq \mathrm{M}_n(F\otimes\HH_{p,\infty})$ for some integer $n$ where $\HH_{p,\infty}$ is the quaternion algebra over $\Q$ ramified only at $p$ and $\infty$. One notes that this is the endomorphism ring of a supersingular elliptic curve in characteristic $p$. 

\subsubsection{} The first section is devoted to the study of $\C$-representations of simple finite dimensional $\Q$-algebras that are rational -- the characteristic polynomials of the elements have coefficients in $\Q$ -- and that admit an invariant symplectic form, we say polarized. We are interested in such representations as both the representation of the endomorphism algebra of an abelian variety given by its Tate module and the representations of $\Q[G]$ coming from the second projection of the map $\iota_{\ell}$ in Definition~\ref{def:pta} have these properties. Our main goal in this section is to show a characterization of such algebras that can be embedded in a nice way into the endomorphism algebras of abelian varieties -- see Definition~\ref{def:goodembed} and Theorem~\ref{theo:goodembed}. In order to do this, we first show in Section~\ref{sub:ratetpolrep} that rational and polarized representations are multiple of a minimal such representation denoted $V_{rp}(E)$. The dimension of this representation depends on the invariants of $E$ and its type in the classification of Albert.  

\medskip

The second section is independent from the first. In the first part of the second section we work on twisted algebras, that is $\Q$-algebras that are obtained by the action of a group on some $\Q$-algebra $E$. Our main interest is to give some restrictions on their centers in the case of the acting group $G$ being cyclic of order $m\geq 1$. Such an algebra is of the form
\[
\bigoplus\limits_{i=0}^{m-1} E.u^i
\]
with, for all $x\in E$, $ux=f(x)u$ for some finite order automorphism $f$ of $E$. When $E$ is simple the center of the twisted algebra is then a product of field extensions of a subfield of $Z(E)$ obtained by extracting roots. A precise statement is given by Proposition~\ref{prop:centrecasspecial1}. We then deal with the case of $E$ not being simple. The second part is dedicated to the proof of Theorem~\ref{theo:centresalgebresintro} which is done by applying the previous results to the case of the rational group algebra of a ramification group $G\simeq \Gamma_p \rtimes \Z/n\Z$. 

\medskip

In the last section we gather together our results and show our main theorem. First we relate the definition of $(p,t,a)$-inertial groups to rational and polarized representations of group algebras. We show that all families of maps $(p_{\ell}\colon G\to \Sp_{2a}(\Q_{\ell}))_{\ell\neq p}$ satisfying the rationality condition of Definition~\ref{def:pta}~$(ii)$ can be recovered from an embedding into the endomorphism algebra of an abelian variety $A_0$ of dimension $a$ over a finite field. We then show that $A_0$ can furthermore be chosen so that the map $G\to \End A_0\otimes \Q$ has values in $\Aut A_0$. We can thus apply Théorème~1.1 of \cite{Ph4} which provides $A$ over a $p$-adic field whose semi-stable reduction $\overline{A}$ is isogenous to $A_0$ and verify that the representation from its $\ell$-adic Tate module recovers $p_{\ell}$ again. We can add the toric part by simply considering $A_0\times \G_m^t$.

\section{Rational representations of simple algebras}

In this section we start by gathering results about representations of finite dimensional simple algebras which are rational and admit a symplectic structure. In the second part of this section we relate these results to algebras coming from the endomorphisms of abelian varieties over finite fields, which we call geometric, and introduce the notion of a good embedding to be used in section \ref{sec:ptaintomonod}. The main result is a characterization of the algebras that admits a good embedding. 

\subsection{Rational and polarizable representations of algebras} \label{sub:ratetpolrep}

\subsubsection{} For the rest of this section let us fix $E$, a simple finite dimensional --since all algebras we are dealing with will be finite dimensional it will be sometimes omitted for clarity -- $\Q$-algebra and denote by $F$ its center. We start by studying the $\C$-representations of such an algebra $E$ and whether those admit invariant bilinear forms. 

\begin{defin}
    A $\C$-representation of $E$ is a $\C$-vector space $V$ equipped with a linear action of $E$, that is a map of $\Q$-algebras 
    \[
    E\longrightarrow \End_{\C} V.
    \]
    For a $\C$-representation $V$ of $E$, let $\operatorname{tr}_V$ denote the central trace of $V$, that is the composition
    \[
    F\to E\to \End_{\C} V \to \C
    \]
    where the last arrow is the trace map. 
\end{defin}

We will also consider $\R$-representations of $E$ in this section, which are defined in the same way.

\begin{lem} \label{lem:classiCrepE}
    The simple $\C$-representations of $E$ are classified by their central trace.
\end{lem}
\begin{proof}
    The simple $\C$-representations of $E$ correspond to simple $E\otimes_{\Q} \C$-modules. For $f$ the degree of $F$ over $\Q$ there is an integer $n$ such that we have an isomorphism
    \[
    E\otimes_{\Q} \C \simeq \prod\limits_{i=1}^f M_n(\C).
    \]
    This isomorphism is such that $F\otimes_{\Q}\C$ is embedded diagonally and by the projection on each factor we recover the embeddings $\sigma\colon F\to \C$ up to multiplication by $n$ by taking the trace.

    \medskip

    As the simple $\C$-representations of $E$ correspond to a choice of such a projection, we see that their central traces are given by the embeddings $F\to \C$ up to multiplication by $n$ and thus they are classified in this way. 
\end{proof}

Let $\{\sigma_1,\dots, \sigma_f\}$ be the set of embeddings $F\to \C$. Then, for $i\in\{1,\dots, f\}$ we will denote by $V_i$ the simple $\C$-representation of $E$ which has a multiple of $\sigma_i$ for its central trace. The $\C$-representation $V_r(E)=\oplus_{i=}^f V_i$ is the minimal non trivial one with rational central trace. 

\begin{defin}
    A positive involution on $E$ is a $\Q$-linear map ${}^*\colon E\to E$ which satisfies, for all $x,y\in E,$
    \[
    (x^*)^*=x,~(xy)^*=y^*x^*
    \]
    and $\operatorname{Tr}_E(xx^*)>0$ for $x\neq 0$. Here, for $x\in E$, the value $\operatorname{Tr}_E(x)$ is the trace of the map given by left multiplication by $x$ in $E$ as a $\Q$-vector space. 
    
    \smallskip

    If $E$ admits a positive involution we say it is polarizable. 
\end{defin}

We give some basic results on positive involutions and the trace map on $E$, which are standard. The fact that an involution is positive can be checked on $E\otimes_{\Q} \R$ considered with the $\R$-linear extension of the involution. For standard results on positive involutions we refer to chapter~8 of \cite{Sch85}. We start with some known facts that are translated in our setting.

\begin{lem}\label{lem:tracescalcul}
    Let $E$ be a simple $\Q$-algebra with center $F$ of degree $f$ over $\Q$. Let $L/F$ be an extension splitting $E$, that is $E\otimes_F L\simeq \mathrm{M}_n(L)$. Let $\operatorname{tr}_L$ be the trace map on $\mathrm{M}_n(L)$ and note that $E\subset \mathrm{M}_n(L)$. Then for all $x\in E$ we have 
    \[
    \operatorname{Tr}_E(x)=\operatorname{tr}(x).
    \]
    Furthermore, considering the embedding $E\subset E\otimes_{\Q} \C$ with $E\otimes_{\Q} \C \simeq \prod_{i=1}^f \mathrm{M}_n(\C)$ with $\operatorname{tr}_{\C}$ the sum of the trace on each factor, then for $x\in E$ we have $\operatorname{tr}_{\C}(x)=f\cdot \operatorname{Tr}_E(x)$.
\end{lem}
\begin{proof}
    The first part is just \cite{Sch85} Lemma~8.5.11. The second part follows as the embedding $E\subset E\otimes_{\Q} \C$ is the diagonal embedding which recovers the embedding $E\subset E\otimes_{F} \C$ on each factor. 
\end{proof}

For a simple $\R$-algebra $\mathrm{M}_n(D)$ let us denote by $\dagger$ the positive involution given by the conjugate-transpose. For a product of simple $\R$-algebra we denote by $\dagger$ again the positive involution given by the conjugate-transpose on each factor.

\begin{lem} \label{lem:transposeconjug}
    Let $E$ be a simple $\Q$-algebra with a positive involution $*$. Then, there is an isomorphism
    \[
    (E\otimes_{\Q} \R, *) \simeq (\prod\limits_{i=1}^r \mathrm{M}_n(D), \dagger). 
    \]
    Furthermore, there is an isomorphism $E\otimes_{\Q} \C\simeq \prod_{i=1}^f \mathrm{M}_n(\C)$ such that $*$ is given by the conjugate-transpose through the corresponding embedding $E\subset \prod_{i=1}^f \mathrm{M}_n(\C)$.
\end{lem}
\begin{proof}
    The proof of the first part is identical to the proof of \cite{Sch85} Theorem 8.13.3. The second part follows from Remark 8.13.4 of \textit{loc. cit.} in the case that $D=\R$ or $D=\HH$. If $D=\C$, the extension of scalars to $\C$ gives a diagonal embedding $\operatorname{M}_n(D)\subset \operatorname{M}_n(\C)\times \operatorname{M}_n(\C)$. It is a simple computation that $*$ is recovered again by the conjugate-transpose on each factor.
\end{proof}

We can now characterize polarizable $\Q$-algebras by the use of hermitian forms. 

\begin{prop}
    The following statements are equivalent.
    \begin{itemize}
        \item[$(i)$] The algebra $E$ is polarizable.
        \item[$(ii)$] All $\C$-representation $V$ of $E$ admits a positive definite hermitian form $H$ such that the adjunction on $\End_{\C} V$ induced by $H$ leaves the image of $E$ invariant.
        \item[$(iii)$] The representation $V_r(E)$ admit a positive definite hermitian form $H$ such that the adjunction on $\End_{\C} V$ induced by $H$ leave $E$ invariant.
    \end{itemize}
\end{prop}

\begin{proof}
    For $(i)\Rightarrow (ii)$, let $V$ be a $\C$-representation of $E$ and consider its decomposition $V=\oplus_{i=1}^f V_i^{n_i}$. The image of $E\otimes_{\Q} \C$ in $\End_{\C} V$ is a diagonal embedding of the simple factors such that $n_i\neq 0$ repeated $n_i$ times. By the choice of a basis of $V$, and thus an isomorphism $\mathrm{M}_n(\C)\simeq \End_{\C} V$ and by Lemma~\ref{lem:transposeconjug} the adjunction to the standard positive definite hermitian form, which is the conjugate-transpose, restricts to the image of $E\otimes_{\Q} \C$ in $\mathrm{M}_n(\C)$ by stabilizing the image of $E$.  

    \medskip

    The implication $(ii)\Rightarrow (iii)$ is trivial. So consider the $\C$-representation $V_r(E)$ and assume it admits a positive definite hermitian form $H$ such that the image of $E$ in $\End_{\C} V_r(E)$ -- note that the sequence of maps $E\to E\otimes_{\Q} \R\to \End_{\C} V_r(E)$ is injective -- is stable by the adjunction induced by $H$. By Lemma~\ref{lem:tracescalcul} this adjunction restricts to a positive involution on $E\otimes_{\Q} \R$ which stabilizes $E$ and thus it induces a positive involution on $E$. 
\end{proof}

For the rest of the section we will furthermore assume that $E$ is polarizable with positive involution ${}^*$. We can now define the dual representation of a $\C$-representation of $E$.

\begin{defin}
    Let $V$ be a $\C$-representation of $E$ given by the action $\rho\colon E\to \End_{\C} V$. Then the dual representation $V^{\vee}$ of $E$ is the dual vector space $V^{\vee}$ with the action of $x\in E$ given by $\rho(x^*)^t$. 

    \smallskip

    We say that a non degenerate bilinear form $B$ on $V$ is invariant if it induces an isomorphism of $\C$-representations $V\simeq V^{\vee}$.

    \smallskip

    We say that $V$ is polarized if there is a non degenerate symplectic form on $V$ which is invariant. 
\end{defin}

Note that for $x\in E$ and $v\in V^{\vee}$, the action of $x$ on $v$ is given, in matrix form with a basis $B$ of $V$ and the corresponding dual basis of $V^{\vee}$ by
\[
\rho(x^*)^t\cdot v = (v^t\cdot \rho(x^*))^t. 
\]

\begin{lem}
    Let $V$ be a $\C$-representation of $E$. There is a positive definite hermitian form on $V$ which induces a $\C$-anti-isomorphism $V\simeq V^{\vee}$. In particular, the central trace $\operatorname{tr}_{V^{\vee}}$ of $V^{\vee}$ is given by $\overline{\operatorname{tr}_V}$.
\end{lem}

\begin{proof}
    As before we consider an hermitian form $H$ on $V$ such that the adjunction coincides with ${}^*$ on the image of $E$ and choose a basis for which it is the conjugate transpose. Let $h\colon V\to V^{\vee}$ be the anti-isomorphism induced by $H$. Then for $x\in E$ and $v\in V$, we have that 
    \[
    h(\rho(x)\cdot v) = \overline{ v}^t\cdot \overline{\rho(x)}^t.
    \]
    Now, since $\rho(x)^t=\rho(x^*)$ we get both claims by writing the product $\overline{v}^t\cdot \overline{\rho(x^*)}$ in the dual basis in a standard way, which corresponds to transposing again.  
\end{proof}

\begin{lem} \label{lem:repdescenteR}
    Let $V$ be a simple $\C$-representation of $E$. Then $\operatorname{tr}_V$ has values in $\R$ if and only if $V\simeq V^{\vee}$ and there is an $\R$-representation $V_0$ of $E$ such that $V\simeq V_0\otimes_{\R} \C$ if and only if there is a non-degenerate bilinear symmetric form on $V$ which is invariant. 
\end{lem}

\begin{proof}
    The first part of the statement follows directly from Lemma~\ref{lem:classiCrepE}. For the second part, if $V\simeq V_0\otimes_{\R} \C$ note that $V_0$ corresponds to a simple factor of $E\otimes_{\Q} \R$ of the form $M_n(\R)$ and not $\mathrm{M}_n(\HH)$ as $V$ would not be simple in that last case. In that case, the positive involution can just be given by the transposition on $\End_{\R} V_0$ and we thus recover a symmetric bilinear form on $V_0$ which induces an isomorphism $V\simeq V^{\vee}$ after base change. For the converse, let $B$ be the invariant bilinear symmetric form on $V$. We follow the proof of \cite{Se77} of the corresponding statement for group representations on p.107--108. Let $H$ be an hermitian positive definite form such that the adjunction recovers the involution on $E$ and $h\colon V\to V^{\vee}$ the induced map. For all $x\in V$ there is a unique element $\varphi(x)$ such that for all $y\in Y$
    \[
    B(x,y)=\overline{H(\varphi(x),y)}.
    \]
    The map $\varphi$ obtained in this way is antilinear and bijective. Its square $\varphi^2$ is an automorphism of $V$. Since $B$ is symmetric we have for all $x,y\in V$
    \[
    H(\varphi^2(x),y)=\overline{B(\varphi(x),y)}=\overline{B(y,\varphi(x))}=H(\varphi(y),\varphi(x)).
    \]
    It follows that $\varphi^2$ is hermitian as $H(\varphi(y),\varphi(x))=\overline{H(\varphi(x),\varphi(y)}$ and thus from the previous equality we get
    \[
    H(\varphi^2(x),y)=\overline{H(\varphi^2(y),x)}= H(x,\varphi^2(y)).
    \]
    From the equality, for all $x\in V$, 
    \[
    H(\varphi(x),\varphi(x))=H(\varphi^2(x),x)
    \]
    we see that $\varphi^2$ is positive definite. We thus get a unique hermitian positive definite $v\in \End_{\C} V$ such that $v^2=\varphi^2$ and $v$ commutes with $\varphi$. Set $\sigma=\varphi\circ v^{-1}$, it is an antilinear endomorphism of $V$. We get $\sigma^2=\operatorname{id}$. Let $V_0$ and $V_1$ be the eigenspaces of $\sigma$ for the values $1$ and $-1$ respectively. Since $\sigma$ is antilinear we have $V_1=iV_0$ and
    \[
    V=V_0\oplus iV_0.
    \]
    But since $B$ and $H$ are invariant for the action of $E$ we get that $V_0$ and $V_1=iV_0$ are stable and $V=V_0\otimes_{\R} \C$ as desired. 
\end{proof}

\begin{prop} \label{prop:simplepolarizedrep}
    Let $V$ be a simple $\C$-representation of $E$ and $n\in \N$. Then exactly one of the three  following possibilities holds.
    \begin{itemize}
        \item[$(i)$] The representation $V$ is not isomorphic to its dual $V^{\vee}$.
        \item[$(ii)$] There is a non degenerate symmetric form on $V$ which is invariant.
        \item[$(iii)$] There is a non degenerate bilinear symplectic form on $V$ which is invariant.
    \end{itemize}

    In the first case, the smallest $\C$-representation of $E$ containing $V^n$ which is polarized is $V^n\oplus (V^{\vee})^n$. In the second case if $V^n$ is polarized then $n$ is even. 
\end{prop}
\begin{proof}
    The distinction between the first case and the two others is just whether $V\simeq V^{\vee}$, or whether $\operatorname{tr}_V$ has values in $\R$. If it does, then there is a non-degenerate bilinear form on $V$ coming from the isomorphism $V\simeq V^{\vee}$. Let $B$ be this form and $B=B_0+B_1$ its decomposition in symmetric and alternate parts. Both induce maps of $\C$-representation between $V$ and $V^{\vee}$ so that by Schur's lemma only one is non zero. 

    \medskip

    In the first two cases we see that $V\oplus V^{\vee}$ is polarized by considering the standard symplectic form on a basis adapted to the direct sum decomposition. Now let $W$ be a representation containing $V^n$ which is polarized. Since $W\simeq W^{\vee}$, in the case that $V$ is not isomorphic to its dual we must have $(V^{\vee})^n\subset W^{\vee}$ but also $V^n$ by the isomorphism $W\simeq W^{\vee}$. 

    \medskip

    In the second case, let $\lambda\colon V^n\simeq (V^{\vee})^n $ be the isomorphism induced from an invariant symplectic form on $V^n$. For a simple sub-representation $W\subset V^n$, $W\simeq V$, we have $\lambda(W)\subset W^{\perp}$ as the symplectic form restricts to the zero form on $W$. It follows that $W\oplus \lambda^{-1}(W^{\vee})$ form a symplectic subspace of $V^n$. By induction we thus get that $n$ is even.  
\end{proof}

\subsubsection{} The skew-fields over $\Q$ of finite dimension which are polarizable have been classified by Albert. We will recall this classification here, see~\cite{Mu08} p.186 (202) for more details. First, we need some results on the classification of simple algebras. By class field theory, skew-fields $E$ over $\Q$ are classified by local invariants $\operatorname{inv}_v E$ at the places $v$ of the center $F$. By the Brauer-Hasse-Noether theorem these local invariants are given by an element of $\Q/\Z$ which is moreover in $\Z/2\Z$ if $v$ is an infinite real place and $0$ if it is a complex place, only finitely many of them are non zero and they sum to zero in $\Q/\Z$. In general we extend this classification by Wedderburn's theorem. That is if $E$ is a simple $\Q$-algebra isomorphic to $\mathrm{M}_n(D)$ for some skew-field $D$ with center $F$, we set $\operatorname{inv}_v(E)=\operatorname{inv}_v(D)$ for a place $v$ of $F$. These invariants classify such algebras $E$ -- see \cite{Sch85} chapter 8 for an introduction to Brauer groups and Brauer equivalence and \cite{Har17} for the computation of Brauer groups of local fields and the Brauer-Hasse-Noether theorem. Now, the Albert classification has four types.  
\begin{itemize}
    \item Type I : The skew-field $D$ is a totally real number field.
    \item Type II : The skew-field $D$ is a quaternion algebra over its center $F$. The number field $F$ is totally real and for every embedding $\sigma\colon F\to \R$ there is an isomorphism
    \[
    \R\otimes_{F} D \simeq  M_2(\R). 
    \]
    \item Type III : The skew-field $D$ is a quaternion algebra over its center $F$. The number field $F$ is totally real and for every embedding $\sigma\colon F\to \R$ there is an isomorphism
    \[
    \R \otimes_{F} D \simeq \HH
    \]
    where $\HH$ is the standard quaternion algebra over $\R$.
    \item Type IV : The center $F$ of $D$ is a totally imaginary quadratic extension of a totally real field $F_0$ with conjugation $\overline{.}$ over $F_0$. The local invariants are such that for any finite place $v$ such that $\overline{v}=v$ then $\operatorname{inv}_v(D)=0$ and otherwise $\operatorname{inv}_v(D)+\operatorname{inv}_{\overline{v}}(D)=0$.
\end{itemize}

By the theorem of Wedderburn again we extend this classification to simple and polarizable $\Q$-algebras of finite dimension. That is if $E\simeq \mathrm{M}_n(D)$ for some skew-field $D$ is polarizable then $D$ is polarizable and the type of $E$ is that of $D$.

\smallskip

We can now restate Proposition~\ref{prop:simplepolarizedrep} in this context.

\begin{prop} \label{prop:albertpol} Let $n\in \N$.
    If $E$ is of type III then every simple representation is polarized. If it is of type IV, for every simple representation $V$ of $E$, the smallest representation containing $V^n$ which is polarized is $V^n\oplus (V^{\vee})^n$. If it is of type I or II, for every simple representation $V$, the representation $V^n$ is polarized if and only if $n$ is divisible by $2$. 
\end{prop}
\begin{proof}
    First we see that if $E$ is of type I, II or III if and only if for every simple representation $V$ the central trace has values in $\R$. This settles the case of type IV with Proposition~\ref{prop:simplepolarizedrep}, as in this case a simple representation $V$ is not isomorphic to its dual.
    
    \medskip

    For type I and II we see that every simple representation $V$ is of the form $V_0\otimes_{\R} \C$ and thus the result follows from Lemma~\ref{lem:repdescenteR} and Proposition~\ref{prop:simplepolarizedrep}.

    \medskip

    If $E$ is of type III we see that the simple modules $V$ can not be written in the form $V\otimes_{\R} V_0$ as we have $\dim_{\C} V=2$ while simple $\R$-representation are of dimension $4$. It follows from Lemma~\ref{lem:repdescenteR} together with Proposition~\ref{prop:simplepolarizedrep} that all such $V$ are polarized. 
\end{proof}

\subsubsection{} We move to consider representations that are rational and classify them. We also classify those that are rational and polarized by the type of the polarizable algebra $E$. 

\begin{defin} \label{defin:reprationelle}
    A $\C$-representation $V$ of $E$ is rational if the characteristic polynomials of the image of the elements of $E$ as endomorphisms of $V$ have coefficients in $\Q$. 
\end{defin}

From Wedderburn's theorem we have that $E\simeq \mathrm{M}_n(D)$ for some skew-field $D/\Q$ of finite dimension and center $F$. Let us set $d^2=[D:F]$ and $f=[F:\Q]$ as well as $\deg(E)=ndf$. 

\medskip

The following is roughly Lemma~1 of chapter II of \cite{ST61}. We record it here for completeness. 

\begin{lem} \label{lem:ratfreecentre}
    Let $V$ be a rational $\C$-representation of $E$. Then $V$ is a free $F\otimes_{\Q} \C$-module. 
\end{lem}
\begin{proof}
    Let again $\{\sigma_1,\dots, \sigma_f\}$ denote the embeddings of $F$ in $\C$. The representation $V$ is given by the direct sum of $V_i^{n_i}$ and for $x\in F$ we have
    \[
    \operatorname{tr}_V(x)=nd\sum\limits_{i=1}^f n_i\sigma_i(x).
    \]
    We thus see that $\operatorname{tr}_V(x)$ is rational for all $x\in F$ if and only if $n_i$ is independent of $i$. That is, we have 
    \[
    V\simeq (\bigoplus_{i=1}^f V_i)^m\simeq (\bigoplus_{i=1}^f \C^{nd})^m
    \]
    as a $F\otimes_{\Q} \C$-module.
\end{proof}

\begin{theo} \label{theo:repratetpol}
  All the rational (resp. and polarized) $\C$-representations of $E$ are multiples of the representation $V_r(E)$ (resp. a unique such representation denoted $V_{rp}(E)$). Furthermore we have that $\dim_{\C} V_r(E)=\deg E$, $V_{rp}(E)\simeq V(E)$ if $E$ is of type III and IV, and $V_{rp}(E)\simeq V_r(E)^2$ if E is of type I and II.
\end{theo}
\begin{proof}
    Let $V$ be a rational $\C$-representation of $E$. By Lemma~\ref{lem:ratfreecentre} we have that $V$ is a free $F\otimes_{\Q} \C$-module, that is
    \[
    V\simeq (\prod\limits_{i=1}^f \C)^m
    \]
    for some $m\in\N$. It follows that $mf=\dim_{\C} V$.
    
    \smallskip

    On the other hand, we have an isomorphism $E\otimes_{\Q} \C\simeq \prod_{i=1}^f E_i$ with $E_i\simeq \mathrm{M}_{nd}(\C)$ a simple algebra whose simple module is $V_i=\C^{nd}$. Hence a decomposition 
    \[
    V\simeq \bigoplus\limits_{i=1}^f V_i^{n_i}.
    \]
    It follows that $V_i^{n_i}$ is isomorphic to $\C^{m}$ as a vector space and thus 
    \[
    ndn_i=m.
    \]
    Finally, we get that $V\simeq V_r(E)^{\frac{m}{nd}}$ as a $\C$-representation of $E$. Note that since $V_r(E)=\bigoplus_{i=1}^f V_i$ and $\dim_{\C} V_i=nd$ the equality $\dim_{\C} V_r(E)=\deg(E)$ holds. 

    \medskip

    For the last part of the statement, first note that from the first part if $V$ is rational and polarized it is a multiple of $V_r(E)$. Now, by Proposition~\ref{prop:albertpol} if $E$ is of type III the equality $V_r(E)=V_{rp}(E)$ holds. If $E$ is of type IV, we have that $f=2f^+$ for some $f^+\in \N$ and that for $i\in \{1,\dots,f\}$ the simple representation $V_i$ is not isomorphic to its dual $V_i^{\vee}$. As we also have that $V_i^{\vee}=V_j$ for some $j\in\{1,\dots,f\}$, up to some renumbering we get
    \[
    V_r(E)=\bigoplus_{i=1}^{f^+} V_i\oplus V_i^{\vee}.
    \]
    Again Proposition~\ref{prop:albertpol} gives that $V_r(E)=V_{rp}(E)$ in this case.
    
    \smallskip

    For the remaining cases, where $E$ is of type I and II, let $V\simeq V_r(E)^m$ be a rational and polarized $\C$-representation of $E$. As the isomorphism $V\simeq V^{\vee}$ from the symplectic form on $V$ respects the action of $E$ it induces an isomorphism between corresponding isotypic components of $V$ and $V^{\vee}$. Since in those cases, every isotypic component is isomorphic to its dual, this gives that the symplectic form on $V$ induces a non degenerate invariant symplectic form on each of its isotypic components $V_i^m$. Proposition~\ref{prop:albertpol} gives that $m$ is divisible by $2$. We thus have that $V\simeq V_r(E)^{2k}$ and $V_{rp}(E)=V_r(E)^2$.    
\end{proof}

For the rest of the text we will keep the notation $V_r(E)$ and $V_{rp}(E)$ for the unique rational or rational and polarized $\C$-representations of a simple polarizable $\Q$-algebra $E$ which generates the others.

\subsection{Geometric Albert algebras and good embeddings} \label{sub:geomalgembed}

\subsubsection{} Let us fix a prime number $p$. We will say that a simple $\Q$-algebra $E$ of finite dimension is quasi-split outside $p$ if $E\otimes_{\Q} \Q_{\ell}$ is quasi-split for all primes $\ell\neq p$ -- an algebra is quasi-split if its a matrix algebra over a field. This is equivalent to the fact that the local invariants of $E$ at all finite places not above $p$ are $0$.

\begin{defin}\label{def:goodembed}
    Let $E$ and $E'$ be polarizable simple $\Q$-algebras. A good embedding is an injective morphism $E\hookrightarrow E'$ such that $V_{rp}(E)=V_{rp}(E')$.
\end{defin}

Note that when there is an embedding $E\hookrightarrow E'$ we have that $\deg(E)\mid \deg(E')$. In the case of a good embedding, from Theorem~\ref{theo:repratetpol} it follows that either $\deg(E)=\deg(E')$ or $2\deg(E)=\deg(E')$. 

\begin{lem} \label{lem:goodembdtype4}
    A polarizable simple algebra $E$ quasi-split outside $p$ admits a good embedding in a polarizable algebra $E'$ quasi-split outside $p$ of type IV if and only if it is not of type III.  
\end{lem}

\begin{proof}
    Consider a good embedding $E\hookrightarrow E'$ and denote by $F$ and $F'$ the centers of $E$ and $E'$ respectively. Let $F_0=F\cap F'$. We get an embedding
    \[
    E\otimes_{F_0} F' \hookrightarrow E'.
    \]
    It follows that $\dim_{\C} V_r(E\otimes_{F_0} F')\leq \dim_{\C} V_{rp}(E')=\deg (E').$ As the embedding $E\hookrightarrow E'$ is good we also have that
    \[
    \epsilon \deg E= \deg E'
    \]
    where $\epsilon$ is $1$ if $E$ is of type III or IV and $2$ if it is of type I or II. Now, if $E$ is of type III, then as $F'$ is a CM field we have an inequality
    \[
    \dim_{\C} V_r(E\otimes_{F_0} F')=\deg(E\otimes_{F_0} F')\geq 2\deg (E)
    \]
    and we get a contradiction
    \[
    \deg E=\deg E'\geq 2\deg E.
    \]

    For the converse, if $E$ is of type I or II we consider the embedding
    \[
    E\hookrightarrow E\otimes_{F} K
    \]
    for a CM field extension $K/F$ of degree $2$. Note that $E\otimes_{F} K$ is a type IV simple algebra quasi-split outside $p$. This embedding is good as
    \[
    V_{rp}(E)=V_r(E)^2=V_{r}(E\otimes_F K)=V_{rp}(E\otimes_F K).
    \]
\end{proof}

\subsubsection{} For a simple abelian variety $A$ over $\overline{\F_p}$ the endomorphism algebra $\End A\otimes \Q$ is a polarizable simple $\Q$-algebra quasi-split outside $p$. By the works of Tate, if $E\simeq \End A\otimes \Q$ for some abelian variety $A$ over $\F_{p^n}$ there is a description of $\operatorname{inv}_v E$ in terms of the valuation at $v$ of the Frobenius of $A$ if $p\mid v$. A precise statement is given by \cite{Tat71} Théorème~1. The algebra $D=\End A\otimes_{\Z} \Q$ is central over $\Q(\pi)$ where $\pi$ is the Frobenius of $A$, it is quasi-split outside $p$ and we have
\[
\begin{cases}
    \operatorname{inv}_v(D)=\frac{v(\pi)}{v(p^n)}\cdot [\Q(\pi)_v \colon \Q_p] \text{ for all places } v\mid p \\
    \operatorname{inv}_v(D)=0 \text{ for all places } v \text{ above a prime } \ell\neq p\\
    \operatorname{inv}_v(D)=\frac{1}{2} \text{ for all real places } v
\end{cases}
\]

Furthermore, the equality $2\dim A= \deg D$ holds. 

\begin{defin}
    A polarizable simple $\Q$-algebra $E$ quasi-split outside $p$ is said to be geometric if there is an abelian variety $A$ over $\overline{\F_p}$ such that
    \[
    E\simeq \End A\otimes \Q.
    \]
\end{defin}

For simplicity, we will say that an algebra $E$ is geometric if it is a polarizable simple $\Q$-algebra quasi-split outside $p$ and it is geometric. 

Let us denote by $\HH_{p,\infty}$ the quaternion algebra over $\Q$ quasi-split outside $p$ such that $\operatorname{inv}_p(\HH_{p,\infty})=\operatorname{inv}_{\infty}(\HH_{p,\infty})=1/2$. The following fact is direct from \cite{Tat71} Exemple~(a) p.97. 

\begin{lem} \label{lem:geotype3supersing}
    Let $E$ be a geometric $\Q$-algebra of type III. Then $E$ is the endomorphism algebra of a power of a supersingular elliptic curve over $\overline{\F_p}$ and we have
    \[
    E\simeq \mathrm{M}_n(\HH_{p,\infty}).
    \]
\end{lem}

\begin{lem} \label{lem:typeIIIreducgoodembed}
    If $E$ is of type III, quasi-split outside $p$ and there is a good embedding $E\hookrightarrow E'$ with $E'$ polarizable, geometric and quasi-split outside $p$ then $E\simeq \mathrm{M}_n(F\otimes \HH_{p,\infty})$ for some totally real number field $F$. 
\end{lem}

\begin{proof}
    By Lemma~\ref{lem:goodembdtype4} we have that $E'$ is of type III. Since $E'$ is geometric, the previous lemma gives an isomorphism
    \[
    E'\simeq \mathrm{M}_{n'}(\HH_{p,\infty}).
    \]
    Let $E\simeq \mathrm{M}_n(D)$ so that $\deg(E)=2nf$ where $f=[F:\Q]$ with $F=Z(E)$ the center of $E$. Since the embedding is good we have $\deg (E)= \deg(E')$ which is
    \[
    2nf=2n'
    \]
    so that $n'=nf$. The centralizer $Z_{E'}(E)$ of $E$ in $E'$ is a simple algebra that contains $F$ and
    from the centralizer theorem, see \cite{Sch85} chapter~8, Theorem~5.4 we have that
    \[
    \dim Z_{E'}(E)\cdot \dim E=d^2f\cdot 4n^2f=\dim E'=4n^2f^2.
    \]
    It follows that $d=1$ and $Z_{E'}(E)=F$. 
    It also follows from \cite{Sch85} chapter~8, Theorem~4.5 that
    \[
    E'\otimes_{\Q} F \sim Z_{E'}(F)=E
    \]
    where $\sim$ denotes the Brauer equivalence -- that the algebras have the same local invariants. This is equivalent to the fact that there is an isomorphism
    \[
    E\simeq \mathrm{M}_n(F\otimes \HH_{p,\infty}).
    \]
\end{proof}

\subsubsection{} In order to prove the main theorem of this section we need to recall some facts about Weil numbers and Honda-Tate theory -- see \cite{Tat71}. 

\begin{defin}
    Let $q=p^n$. A Weil $q$-number (sometimes $q$-Weil number) is an algebraic integer $\pi$ such that for all embedding $\sigma\colon\Q(\pi)\hookrightarrow \C$ we have
    \[
    |\sigma(\pi)|=q^{\frac{1}{2}}.
    \]
\end{defin}

The theorem of Honda-Tate states that there is a one-to-one correspondence between conjugacy classes of Weil $q$-numbers and simple abelian varieties over $\F_q$ up to isogeny. In particular, to construct abelian varieties over finite fields it is enough to construct such algebraic integers. 

\begin{prop} \label{prop:type4geomembed} A simple polarizable $\Q$-algebra of type IV quasi-split outside $p$ admits a good embedding in a geometric $\Q$-algebra.
\end{prop}
\begin{proof}
    Let $E$ be a simple polarizable $\Q$-algebra of type IV quasi-split outside $p$. As before, let $E\simeq \mathrm{M}_n(D)$ for some skew-field $D$ over $\Q$ with center $F$. The skew-field $D$ is itself a simple polarizable $\Q$-algebra of type IV with center $F$ and if it admits a good embedding in a geometric polarizable $\Q$-algebra $E'$ of type IV quasi-split outside $p$ then $E\hookrightarrow \mathrm{M}_n(E')$ is a good embedding. So we will assume $E=D$.
    
    \medskip
    
    From the definition the fact that $D$ is of type IV we have that $\operatorname{inv}_v(D)+\operatorname{inv}_{\overline{v}}(D)=0$ in $\Q/\Z$ for all finite places $v$ of $F$. For all those places $v$ above $p$, let us choose a non-negative representative in $\Q$ of $\operatorname{inv}_v(D)$ such that $\operatorname{inv}_v(D)+\operatorname{inv}_{\overline{v}}(D)=[F_v:\Q_p]$. Note that if $v=\overline{v}$ we have that $[F_v:\Q_p]$ is even by the characterization of type IV so that we can avoid choosing $\operatorname{inv}_v(D)=1/2$.  For a place $v$ of $F$ let $e_v$ be the ramification index at $v$. Let $n\in \N\setminus\{0\}$ be such that for all places $v$ above $p$ 
    \[
    n\cdot \frac{\operatorname{inv}_v(D)}{[F_v:\Q_p]}\in \N
    \]
    and for a place $v$ above $p$ let 
    \[
    n_v=ne_v\cdot \frac{\operatorname{inv}_v(D)}{[F_v:\Q_p]}.
    \]
    Note that we have $n_v+n_{\overline{v}}=ne_v$ by construction.

    \smallskip
    
    Consider the ideal $\prod_{v\mid p} \mathfrak{p}_v^{n_v}$ where $\mathfrak{p}_v$ is the prime of $\OO_F$ associated to the place $v\mid p$. Then, up to multiplying $n$ by an integer big enough we can assume this ideal is principal, given by $(\pi)$. We have
    \[
    (\pi\overline{\pi})=\prod\limits_{v\mid p} \mathfrak{p}_v^{n_v+n_{\overline{v}}}=\prod\limits_{v \mid p} (\mathfrak{p}_v^{e_v})^n= (p^n).
    \]
    Thus we have $\pi\overline{\pi}=p^nu_0$ with $u_0\in \OO_{F^+}^{\times}$. Let $m\in \N$ be such that $u_0^m=u\overline{u}$ with $u\in \OO_F^{\times}$, this is possible as the norm map has finite index. 

    \smallskip

    Now consider $\pi_1=\pi^mu^{-1}$. We have 
    \[
    \pi_1\overline{\pi_1}=p^{nm}u_0^{m}u_0^{-m}=p^{mn}.
    \]
    In particular, $\pi_1$ is a Weil $p^{mn}$-integer and the associated abelian variety by Honda-Tate theory $A_0$ is such that $D_0=\End A_0\otimes \Q$ has center $F_0=\Q(\pi_1)$ and invariants 
    \[ \operatorname{inv}_{v_0}(D_0)=\frac{v_0(\pi_1)}{v_0(p^{mn})}\cdot [{F_0}_{v_0}:\Q_p]=\frac{mn_v}{mne_v} \cdot [{F_0}_{v_0}:\Q_p]=\frac{\operatorname{inv}_v(D)}{[F_v:{F_0}_{v_0}]}
    \]
    for all places $p\mid v_0 \mid v$ of $\Q(\pi_1)$ above $p$ and below $v$. 

    \medskip

    It follows that $D=D_0\otimes_{\Q(\pi_1)} F$ and that $D\subset \End A_0^{[F:\Q(\pi_1)]}\otimes \Q=\mathrm{M}_{[F:\Q(\pi_1)]}(D_0)$. The simple algebra $D_0$ has degree $\deg(D_0)=d\cdot [\Q(\pi_1):\Q]$ so that
    \[
    \deg(D)=df=d\cdot [F:\Q(\pi_1)]\cdot  [\Q(\pi_1):\Q]=\deg(\mathrm{M}_{[F:\Q(\pi_1)]}(D_0))
    \]
    and the embedding is good.
\end{proof}

We can now state the main theorem of this section. 

\begin{theo} \label{theo:goodembed}
    A simple and polarizable $\Q$-algebra $E$ quasi-split outside $p$ admits a good embedding into a geometric $\Q$-algebra $E'$ if and only if it is either not of type III or of the form $\mathrm{M}_n(F\otimes \HH_{p,\infty})$ for some totally real number field $F$.
\end{theo}

\begin{proof}
    Assume first that $E$ is of type III. Then by Lemma~\ref{lem:typeIIIreducgoodembed} if $E$ admits a good embedding it is of the required form. So assume furthermore that $E\simeq \mathrm{M}_n(F\otimes \HH_{p,\infty})$ for some totally real number field $F$. Then the embedding
    \[
    E\hookrightarrow \mathrm{M}_{fn}(\HH_{p,\infty})
    \]
    is good where $f=[F:\Q]$.
    
    \medskip

    In all other cases, combine Lemma~\ref{lem:goodembdtype4} with Proposition~\ref{prop:type4geomembed} to get the result.
\end{proof}

\section{The rational group algebra of ramification groups}

We start by developing some technical points on $\Q$-algebras that are twisted by cyclic groups. The main results is given in the second part of the section and applies the results of the first part to give a description of the rational group algebra $\Q[G]$ of ramification groups, giving precisions to a theorem of Serre and generalizing results of Roquette and Ford.  

\subsection{Twisted $\Q$-algebras} 

\subsubsection{} \label{subsub:twisted1} Let $G$ be a finite group. For a finite dimensional $\Q$-algebra $E$ and a morphism $\varphi\colon G\to \Aut E$ we denote by $E[G,\varphi]$ the $\Q$-algebra defined by $\bigoplus_{g\in G} E.g$ where the multiplication is given, for all $x,y\in E$ and all $g,h\in G$, by
\[
xg\cdot yh= x\varphi(g)(y)\cdot gh.
\]
When $E$ has a positive involution $*$ and $\varphi$ verifies $\varphi(g)(x)^*=\varphi(g)(x^*)$ for all $g\in G$ and $x\in E$, that is $\varphi(G)\subset \Aut (E,*)$, then there is an induced involution on $E[G,\varphi]$ given by 
\[
(xg)^*=g^{-1}x^*=\varphi(g^{-1})(x^*)g^{-1}.
\]

\begin{lem}
    The involution induced by $*$ on $E[G,\varphi]$ is positive if and only  if $*$ is.
\end{lem}
\begin{proof}
    For all $x\in E$ and $g\in G$ we have
    \begin{align*}
        (xg)(xg)^*&=xg \varphi(g^{-1})(x^*)g^{-1} \\
        &=x\varphi(g)\big(\varphi(g^{-1})(x^*)\big )gg^{-1}\\
        &= xx^*.
    \end{align*}
    The statement now follows as the multiplication by $xx^*$ on $E[G,\varphi]$ stabilizes the direct sum decomposition we have $\operatorname{Tr}_{E[G,\varphi]}(xx^*)=\Card G\cdot \operatorname{Tr}_E(xx^*)$.
\end{proof}

Note that for a finite group $H$ the rational group algebra $\Q[H]$ is equipped with a natural positive involution -- see the last section of \cite{Dav06} for example -- given for $h\in H$ by
\[
h\mapsto h^{-1}.
\]
Thus from Proposition~8.1 of \cite{Ré17} each of the simple factor $E\subset \Q[G]$ that appear in the decomposition of $\Q[G]$ is a polarizable $\Q$-algebra. 
This is a particular case of the construction presented here with $\varphi\colon G\to \Aut \Q$ the trivial map and $*$ the identity on $\Q$.

\smallskip

Assume now that $H$ is carrying a $G$ action by a map $G\to \Aut H$. The action extends to a map $\varphi\colon G\to \Aut \Q[H]$ and we have 
\[
\Q[H][G,\varphi]= \Q[H\rtimes G].
\]
One can check that the action of $G$ is compatible with the natural involution on $\Q[H]$ and  the induced involution on $\Q[H\rtimes G]$ is again the natural one.

\subsubsection{} We will now only consider the case when the acting group $G$ is cyclic. Let $m\in \N$ be such that $G\simeq \Z/m\Z$ and let $\sigma=\varphi(1)\in \Aut E$. In this situation we will write $E[m,\sigma]$ in place of $E[G,\varphi]$. We also write 
\[
E[m,\sigma]= \bigoplus\limits_{i=0}^{m-1} E.u^i
\]
with $ux=\sigma(x)u$ for $x\in E$ and $u^m=1$. 

\medskip

We are interested in the center of such algebras. Let us start by a special case. The algebra $E$ is simple and let $\tau=\sigma_{|Z(E)}\in \Aut Z(E)$. The order $n$ of $\tau$ divides $m$ and we have that $\sigma^n$ is an automorphism of the simple algebra $E$ central over $F=Z(E)$. By Skolem-Noether's theorem we get that there is an element $y\in E^{\times}$ such that for all $x\in E$ we have 
\[
\sigma^n(x)=y^{-1} xy.
\]
We show that we can further choose $y$ to be fixed by $\sigma$ and have further properties.

\begin{prop} \label{prop:centrecasspecial1}
    Let $E$ be a simple algebra with an automorphism $\sigma$ of order $m$. Let $\tau=\sigma_{\mid Z(E)}$ and $n$ its order. Then we can choose $y\in E^{\times}$, such that for all $x\in E$
    \[
    \sigma^n(x)=y^{-1}xy,
    \]
    to be fixed by $\sigma$. In this case, $y^{\frac{m}{n}}$ is an element of $F^{\tau}$ and there is an isomorphism
    \[
    \psi\colon Z(E[m,\sigma]) \simeq F^{\tau}[X]/(X^{\frac{m}{n}}-y^{\frac{m}{n}}).
    \]
    Moreover, if $E$ has a positive involution $*$ stable by $\sigma$ which induces an involution noted $*$ again on $E[m,\sigma]$ we can choose $\psi$ such that $X^*=yy^*X^{-1}$ in the quotient. In particular, we have $yy^*\in F^{\tau}$.
\end{prop}
\begin{proof}
    Let us start by showing we can choose such a $y$ to be fixed by $\sigma$. Start with $y\in E^{\times}$ given by Skolem-Noether's theorem applied to $\sigma^n$. As the automorphisms $\sigma$ and $\sigma^n$ of $E$ commute, we have for all $x\in E$
    \[
    \sigma(\sigma^n(x))=\sigma(y)^{-1}\sigma(x) \sigma(y)=\sigma^n(\sigma(x))=y^{-1}\sigma(x) y
    \]
    so that $\sigma(x)\sigma(y)y^{-1}=\sigma(y)y^{-1} \sigma(x)$ and thus $t=\sigma(y)y^{-1}$ is an element of $F$. 

    \smallskip

    In the cyclic extension $F/F^{\tau}$ the norm of $t$ is 
    \[
    N_{F/F^{\tau}} (t) = \prod\limits_{i=0}^{n-1} \tau^i(t)=\prod\limits_{i=0}^{n-1} \sigma^i(\sigma(y)y^{-1})=\sigma^n(y)y^{-1}=1.
    \]
    By the Hilbert's theorem 90 we get an element $s\in F^{\times}$ such that $t=\tau(s)s^{-1}$. Let us set $y'=ys^{-1}$. For all $x\in E$ we have
    \[
    \sigma^n(x)={y'}^{-1}xy'
    \]
    as $s$ commutes with $x$. The element $y'$ verifies further
    \[
    \sigma(y')=\sigma(y)\sigma(s^{-1})=\sigma(y)\sigma(t)\sigma^2(s)^{-1}=\sigma^2(y)\sigma^2(s)^{-1}=\sigma^2(y')
    \]
    which is just $y'=\sigma(y')$ as desired. For all $x\in E$ we also have
    \[
    {y'}^{-\frac{m}{n}} x {y'}^{\frac{m}{n}}= \sigma^{n\cdot \frac{m}{n}}(x)=x.
    \]
    In particular, ${y'}^{\frac{m}{n}}$ is central and fixed by $\tau$. 

    \medskip

    We will now assume that $y$ was chosen such that it has the properties proven in the first part and we show that there is an isomorphism $\psi\colon E[m,\sigma] \simeq F^{\tau}[X]/X^{\frac{m}{n}}-y^{\frac{m}{n}}$. By definition we have
    \[
    E[m,\sigma]= \bigoplus\limits_{i=0}^{m-1} Eu^i.
    \]
    Let us consider an element $z\in Z(E[m,\sigma])$ given by 
    \[
    z= \sum\limits_{i=0}^{m-1} a_i u^i.
    \]
    Since $zu=uz$ we get $\sigma(a_i)=a_i$ for all $i\in \{0,\dots, m-1\}$ and the equality $zx=xz$ for all $x\in E$ gives $a_i \sigma^i(x)=xa_i$. For $i\in \{0,\dots, m-1\}$ such that $n$ does not divide $i$ there is an $x\in F$ such that $\sigma^i(x)\neq x$ and thus $a_i=0$. Let us set $b_i=a_{ni}y^{-i}$ for $0\leq i \leq \frac{m}{n}-1$. For all $x\in E$ we have 
    \[
    xb_i=xa_{ni}y^{-i}=a_{ni}\sigma^{ni}(x)y^{-i}=a_{ni}y^{-i}x=b_ix
    \]
    so that $b_i\in F$ and $\sigma(b_i)=b_i$ gives furthermore that $b_i\in F^{\tau}$. Going back to the element $z$ we have
    \[
    z=\sum\limits_{i=0}^{\frac{m}{n}-1} b_iy^i u^{ni}
    \]
    with $b_i\in F^{\tau}$. Conversely let us show that all elements $z$ of this form are in the center of $E[m,\sigma]$. It is enough to check that it is the case for $yu^n$. Let $x\in E$, we have
    \[
    yu^nx=y\sigma^n(x)u^n=xyu^n
    \]
    and we also have
    \[
    uyu^n=\sigma(y)u^{n+1}=yu^n\cdot u.
    \]
    Thus we have $Z(E[m,\sigma])=F^{\tau}[yu^n]$ with $(yu^n)^{\frac{m}{n}}=y^{\frac{m}{n}}$. 

    \medskip

    For the last part, note first that for all $x\in E$ 
    \[
    ({yy^*})^{-1}xyy^{*}={y^*}^{-1}y^{-1} x y y^{*}={y^{*}}^{-1} \sigma^n(x)y^*=(y\sigma^n(x^*)y^{-1})^*=x
    \]
    and thus $yy^*\in F^{\tau}$. We now compute 
    \[
     {yu^n}{(yu^n)}^*= yu^nu^{-n} y^*=yy^*.
    \]
\end{proof}

Specifying the context again we get a more precise statement.

\begin{prop} \label{prop:centralspecial2}
    Assume that $E\simeq \mathrm{M}_d(K)$ where $K$ is a CM field and $*$ is a positive involution on $E$ stable by $\sigma$. Let us assume  furthermore that $\tau$ is the complex conjugation on $K$ and that we have chosen $y$ as in Proposition~\ref{prop:centrecasspecial1}. Then $yy^*\in F^{\tau}$ is the square of an  element  of $K$.
\end{prop}

\begin{proof}
    Let us denote by $\overline{\cdot}$ the complex conjugation and its extension to $E$ acting on the coefficients on the matrices. By conjugation by an element $w\in E^{\times}$ we have $*$ given for all $x\in E$ by $x^*=w^{-1}\overline{x}^tw$. By Skolem-Noether's theorem we also have that $\sigma(x)=z^{-1}\overline{x}z$ for an element $z\in E^{\times}$ and all $x\in E$. 

    \medskip

    As for $x\in E$ we have
    \begin{itemize}
        \item[$(i)$] $(x^*)^*=x$, which is $(w^{-1}\overline{x}^t w)^*=w^{-1} \overline{w}^tx(w^{-1} \overline{w}^t)^{-1}=x$
        \item[$(ii)$] $\sigma(x)^*=\sigma(x^*)$ which is $w^{-1}\overline{\sigma(x)}^tw=w^{-1}\overline{z}^tx \overline{z^{-1}}^tw=z^{-1}\overline{w^{-1}} x^t \overline{w}z$ 
        \item[$(iii)$] $\sigma^2(x)=y^{-1}xy$ which is $\sigma(z^{-1}\overline{x}z)=z^{-1}\overline{z^{-1}}x\overline{z}z=y^{-1}xy$
    \end{itemize}
    the elements $w^{-1}\overline{w}^t$, $\overline{w}zw^{-1} \overline{z^{-1}}$ and $\overline{z}zy^{-1}$ are central. By $\sigma(y)=y$ we even have that $y^{-1}z\overline{z}$, which is an element of $K$, is fixed by the complex conjugation and thus by $*$. We now have
   \begin{align*}
    yy^*&=(\overline{z}zy^{-1})^{-1}\cdot \overline{z}z  \cdot (\overline{z}z \cdot \overline{z}zy^{-1})^{-1})^*=(\overline{z}zy^{-1})^{-2}\cdot \overline{z}z\cdot (\overline{z}z)^* \\
    &= (\overline{z}zy^{-1})^{-2}\cdot \overline{z}z \cdot w^{-1} \overline{z}^t w \cdot w^{-1} z^t w=(y^{-1} z\overline{z})^{-2}\cdot \overline{z}z \cdot w^{-1} \overline{z}^t z^t w \\
    &= (\overline{z}zy^{-1})^{-2}\cdot \overline{z} \overline{w^{-1}}\overline{w}zw^{-1} \overline{z}^t z^t w=(\overline{z}zy^{-1})^{-2} \cdot \overline{w}zw^{-1}\overline{z}^t \cdot \overline{z} \overline{w^{-1}} z^tw\\
    &=(\overline{z}zy^{-1})^{-2}\cdot \overline{w}zw^{-1}\overline{z}^t \cdot \overline{z}\overline{w^{-1}} w^t {w^{-1}}^tz^t \overline{w}^t \overline{w^{-1}}^tw. \\
   \end{align*}
   Notice now that since $w^{-1}\overline{w}^t$ is central we have $w^{-1}\overline{w}^t=\overline{w}{w^{-1}}^t$ and thus
   \begin{align*}
       yy^*&=(\overline{z}zy^{-1})^{-2}\cdot \overline{w}zw^{-1}\overline{z}^t \cdot \overline{z}{w^{-1}}^tz^t\overline{w}^t \cdot (w^{-1} \overline{w}^t)^{-2}\\
       &=(\overline{z}zy^{-1})^{-2}\cdot (\overline{w}zw^{-1}\overline{z}^t)^2 \cdot{\overline{w}^t}^{-2}.
   \end{align*}
   The last equality is again from the fact that $\overline{w}zw^{-1}\overline{z}^t$ is central and thus fixed by the transposition. We finally have obtained that $yy^*$ is a square of an element of $K$. 
\end{proof}

\subsubsection{} In order to deal with the general case, where $E$ is not assumed simple, we will need a few lemmas that might be well known.

\begin{lem} \label{lem:decompactionalg}
    Let $E$ be a simple $\Q$-algebra, $k,m\geq 1$ and $f\in \Aut E^k$ such that $f^m=1$. Then we can choose a basis $(x_1,\dots,x_k)$ and a partition $I_1,\dots, I_r$ of $\{1,\dots, k\}$ such that for every $j\in \{1,\dots, r\}$ the action of $f$ stabilizes $\bigoplus_{i\in I_j} Ex_i$ and acts on such a factor by 
    \[
    f(x_1,\dots, x_{\Card I_j})=(\sigma_j(x_{\Card I_j}), x_1,\dots, x_{\Card I_j-1})
    \]
    where $\sigma_j\in \Aut E$ and we have chosen a renumbering of the $x_i's$ such that $i\in I_j$. 
\end{lem}

\begin{proof}
 First note that $f$ acts on $\Q^k\subset E^k$ by permutation. We consider the partition $I_1,\dots, I_r$ of $\{ 1,\dots, k\}$ obtained from the decomposition of the permutation action of $f$ on $\Q^k$ as disjoint cycles. For $j\in \{1,\dots, r\}$ we can now consider the algebra -- which is a subset of $E^k$ -- which corresponds to the subset $I_j$, is isomorphic to $E^{\Card I_j}$ and which is stable by $f$. In this subalgebra the images of $E\times 0\cdots \times 0$ are all isomorphic to $E$ and we have an isomorphism
 \[
 E^{\Card I_j} \simeq E\times f(E)\times \cdots \times f^{\Card I_j-1}(E).
 \]
 In a basis adapted to this isomorphism we get the result. As this holds for all $j\in \{1,\dots, r\}$, we are done.

\end{proof}

\begin{lem}
    Let $E$ be a $\Q$-algebra, $m,k\geq 1$ and $\sigma\in \Aut E$ of order dividing $m$. Then we have an isomorphism
    \[
    E[m,\sigma]^k= E^k[m, \sigma^{\times k}]
    \]
    where $\sigma^{\times k}\in \Aut E^k$ is the product. 
\end{lem}

\begin{proof}
    By definition we have
    \[
    E[m,\sigma]^k= \bigoplus\limits_{j=1}^k \bigoplus\limits_{i=0}^{m-1}E.u^i=\bigoplus\limits_{i=0}^{m-1} \bigoplus\limits_{j=1}^k E.u^i.
    \]
    Now, since $\bigoplus_{j=1}^k E.u^i$ is just $E^k(u^{\times k})^i$ where $u^{\times k}$ is the $k$-tuple with each coordinate being $u$ we thus have
    \[
    E[m,\sigma]^k = E^k[m,\sigma^{\times k}].
    \]
\end{proof}

We go back to the situation of Lemma~\ref{lem:decompactionalg} to show that we are left with matrix algebras over some algebra of the form $E[m,\sigma]$ for each stable factor. 

\begin{lem} \label{lem:decompalgebrematrices}
    Let $E$ be a simple algebra, $m,k\geq 1$ and $f\in \Aut E^k$. Let us assume that $f$ acts on $E^k$ by
    \[
    f(x_1,\dots, x_k)= (\sigma(x_k),x_1,\dots ,x_{k-1})
    \]
    for some $\sigma\in \Aut E$ with $\sigma^m=1$. Then we have
    \[
    E^k[mk,f]\simeq \mathrm{M}_{k}(E[m,\sigma]).
    \]
\end{lem}

\begin{proof}
    Let $\Delta$ be the subalgebra of diagonal matrices in $\mathrm{M}_k(E[m,\sigma])$. By definition we have
    \[
    \Delta \simeq E[m,\sigma]^k.
    \]
    
    Now consider the matrix of $M_k(E[m,\sigma])$ defined as
    \[
    M=\begin{bmatrix}
        0 &  &  & u \\
        1 & \ddots  &  &  \\
         & \ddots & \ddots & \\
         &  & 1 & 0
    \end{bmatrix}
    \]

    We have

    \[
    \mathrm{M}_k(E[m,\sigma)]=\bigoplus\limits_{i=0}^{k-1} \Delta M^i= \bigoplus\limits_{i=0}^{k-1} E[m,\sigma]^k M^i.
    \]

    By the previous lemma, we have $E[m,\sigma]^k\simeq E^k[m,\sigma^{\times k}]$. Hence we have
    \begin{align*}
        \mathrm{M}_k(E[m,\sigma)]& = \bigoplus\limits_{i=0}^{k-1} E^k[m,\sigma^{\times k}] M^i \\
        &=\bigoplus\limits_{i=0}^{k-1} \bigoplus\limits_{j=0}^{m-1} E^k u^j M^i
    \end{align*}
    where on the last equality $u$ is identified as the diagonal matrix $uI$. But since $uI=M^k$ we have
    \[
    \mathrm{M}_k(E[m,\sigma)]=\bigoplus\limits_{i=0}^{k-1}\bigoplus\limits_{j=0}^{m-1} E^k  M^{kj+i}.
    \]
    Let us now remark that by construction, for $x\in E^k$, we have $Mx=f(x)M$ so that
    \[
    \mathrm{M}_k(E[m,\sigma)]=\bigoplus\limits_{i=0}^{k-1}\bigoplus\limits_{j=0}^{m-1} E^k  M^{kj+i}=\bigoplus\limits_{i=0}^{km-1} E^kM^i=E^k[mk,f].
    \]
\end{proof}

Our last lemma concerns the general case but with an interior automorphism.

\begin{lem} \label{lem:isomcentrealgebreint}
    Let $E$ be a $\Q$-algebra, $m\geq 1$ and $f\in \operatorname{Int}E$ such that $f^m=1$. Then we have an isomorphism
    \[
    E[m,f] \simeq E\otimes_{Z(E)} Z(E[m,f]).
    \]
\end{lem}
\begin{proof}
 Since $f$ is interior, let $y\in E^{\times}$ be such that for all $x\in E$
 \[
 f(x)=y^{-1}xy.
 \]
 As $f^m=1$ we have $y^m\in Z(E)$. We further have that in $E[m,f]$ the element $yu$ is central. Indeed, for all $x\in E$ we have
 \[
 xyu=yf(x)u=yf(x)y^{-1}yu=xyu
 \]
 and
 \[
 uyu= f(y)u^2=yu\cdot u.
 \]
 We thus have
 \[
 E[m,f]= \bigoplus\limits_{i=0}^{m-1} E(yu)^i=E\otimes_{Z(E)} \bigoplus\limits_{i=0}^{m-1} Z(E)(yu)^i.
 \]
 From this we get 
 \[Z(E[m,f])= \bigoplus\limits_{i=0}^{k-1} Z(E)(yu)^i\]
 and the result follows. 
\end{proof}

We can now show that the centers of the $\Q$-algebras we are interested in have a specific form.

\begin{theo} \label{theo:centresalgebres}
    Let $E$ be a simple algebra with positive involution $*$ and $m,k\geq 1$. Let $f\in \Aut (E^k, *^{\times k})$ with $f^m=1$. Assume that $Z(E)\subset K$ for some CM field $K$ and that one of the following holds
    \begin{enumerate}
        \item $m$ is odd
        \item $E\simeq \mathrm{M}_d(K)$
        \item $E\simeq \mathrm{M}_d(\Q)$.
    \end{enumerate}

    Then $Z(E^k[m,f])$ is a product of CM fields and subfields of $K$ or in case (3) a product of such fields and of quadratic extensions of $\Q$. Furthermore, in case (3), these are matrix algebras over their centers.
\end{theo}
\begin{proof}
    From Lemma~\ref{lem:decompactionalg} the action breaks into disjoint part over which Lemma~\ref{lem:decompalgebrematrices} gives a description in terms of matrix algebras of the form $\mathrm{M}_r(E[m',\sigma])$ for some $m'\mid m$, $r\leq k$ and $\sigma\in \Aut E$. The center of such an algebra is given by the center of $E[m',\sigma]$ so we need only to deal with the case $k=1$.

    \medskip

    Let us assume $k=1$. We are thus in the situation of Proposition~\ref{prop:centrecasspecial1} whose notation we reuse. Each factor of the center is thus an extension $L/F^{\tau}$ generated by an element $\alpha$ such that $\alpha^{\frac{m}{n}}\in F^{\tau}$ and $\alpha\alpha^*=yy^*$. Since $E[m,\sigma]$ is equipped with a positive involution, each factor of its center, and $L$ in particular, is also equipped with such an involution. By the classification of Albert they are thus CM or totally real.

    \smallskip
    
    If $L$ is not CM, it is totally real and so is $F^{\tau}$. In this case it follows that $\alpha\alpha^*=\alpha^2\in F^{\tau}$. Now, if (1) holds, note that $m/n$ is odd so that $\alpha^{\frac{m}{n}}\in F^{\tau}$ and $\alpha^2\in F^{\tau}$ give $\alpha\in F^{\tau}$ which is just $L=F^{\tau}$.

    \medskip

    In the case of (2), we have $F^{\tau}=K^{\tau}$ a totally real subfield of $K$ so that $\tau^{\frac{n}{2}}$ is the complex conjugation. Applying Proposition~\ref{prop:centralspecial2} with $\sigma=f^{\frac{n}{2}}$ we get $\beta\in K$ such that $\alpha^2=yy^*=\beta^2$ and $L=K^{\tau}(\beta)\subset K$. 

    \medskip

    In the last case, we just note again that if $L$ is not CM it is a totally real field, quadratic over $\Q$. The last statement follows from Lemma~\ref{lem:isomcentrealgebreint} as in this case the action for $k=1$ must be interior.
\end{proof}

\subsection{The rational group algebra of ramification groups}

\subsubsection{} Let $p$ be a prime number. A ramification group at $p$ is a finite group of the form $\Gamma_p\rtimes \Z/n\Z$ where $\Gamma_p$ is a $p$-group and $n$ is prime to $p$. We will omit to say "at $p$" when the prime $p$ is clear from context. It is part of the main result of \cite{Mau68} that all such groups can be realized as the inertia group of an extension of $p$-adic fields. Recall from paragraph~\ref{subsub:twisted1} that the group algebra $\Q[G]$ and each of its simple factors are equipped with a natural positive involution. 

\medskip

We start by a lemma in the specific case where $\Card G$ is odd and it verifies condition $(B_q)$ of \cite{Se60} p.409 for an odd prime $q\neq p$. 
\[
(B_q)~\textit{There exists a normal subgroup $C$ of $G$ such that } G/C \textit{ is a } q{\textit{-group}}.
\]

\begin{lem} \label{lem:notypeIIIoddorder}
    Let $G$ be a ramification group of odd order and satisfying $(B_q)$ for an odd prime $q\neq p$. Then the algebra $\Q[G]$ has no simple factor of type III.
\end{lem}
\begin{proof}
    We follow the description given by the second case in the proof of Théorème~3 of \cite{Se60}. With the fact that $G$ is a ramification group at $p$, satisfies $B_q$ and is of odd order we have an isomorphism
    \[
    G\simeq \Z/p^e\Z \rtimes (\Z/m\Z\times \Z/q^f\Z)
    \]
    for some integers $e,f,m\geq 0$ and $m$ prime to $p$ and $q$. Let $T$ be the kernel of the conjugation action of $G$ on $\Z/p^e\Z$. By definition, $\Z/p^e\Z\subset T$ but also $C\subset T$ since it is a cyclic group containing $\Z/p^e\Z$. We get that $G$ is an extension of $T\simeq \Z/p^e\Z \times \Z/mq^w\Z$, $w\geq 1$ an integer less than $f$, and $R=G/T\simeq \Z/q^{f-w}\Z$. This extension is given by an element $a\in H^2(R,T)$ with values in $\Z/mq^w\Z \subset T$. 

    \medskip

    Now, we have $\Q[G]=\Q[T][R,\varphi]$ for the action $\varphi\colon R\to \Aut T$. As $\Aut T$ stabilizes the decomposition of $\Q[T]$ into a product of cyclotomic fields it is also the case of the action of $R$. We get 
    \[
    \Q[G]\simeq \bigoplus\limits_{d\mid p^emq^w} A_d
    \]
    where $A_d$ is $\Q(\mu_d)[R,\varphi']$. For $d= \Card T$ we denote $A_d$ by $A[G]$. The algebra $A[G]$ is a cyclic algebra given by $(\Q(\mu_{p^emq^w})/F, a)$. All other simple factors are described in the same way as $A[G/G']$ for a normal subgroup $G'\subset G$ and thus satisfies the conditions of the lemma.

    \medskip

    By construction, the field $F\subset \Q(\mu_{p^emq^w})$ is the center of $A[G]$ and contains $\Q(\mu_{mq^w})$ as the action of $R$ is trivial on the corresponding subgroup of $G$. If $F$ is totally real we must have $mq^w\in \{1,2\}$. Since the order of $G$ is odd this implies $mq^w=1$ but then $a$ which has values in $\Z/mq^w\Z$ is just trivial, hence $A[G]$ is of type II. As this applies to all simple factors of $\Q[G]$, from the fact they are all of the form $A[G/G']$, we have that $\Q[G]$ has no factor of type III.  
\end{proof}

\subsubsection{} We can now prove our main technical result which is the following structure theorem. It gives precision to Théorème~3 of \cite{Se60} on the group algebra $\Q[G]$ and generalizes the results of Roquette and Ford on $p$-groups in \cite{Roq58, Fo87}.

\begin{theo} \label{theo:structureQG}
    Let $G$ be a ramification group at $p$. The rational group algebra $\Q[G]$ has the following properties. 
    \begin{itemize}
        \item[$(i)$] The algebra $\Q[G]$ is quasi-split outside $p$.
        \item[$(ii)$] If $E\subset \Q[G]$ is a simple factor of $\Q[G]$ then the center $Z(E)$ of $E$ is a CM field or a subfield of $\Q(\mu_{p^{\infty}})$.
    \end{itemize}
    Furthermore, if $\Card G$ is odd then no simple factor of $\Q[G]$ is of type III.  
\end{theo}
\begin{proof}
    The statement $(i)$ is exactly Théorème~3 of \cite{Se60} by the equivalence given in Proposition~1 of \textit{loc. cit}.

    \medskip

    Let us prove $(ii)$. Since $G$ is a ramification group it is of the form $H\rtimes \Z/m\Z$ with $H$ a $p$-group and $m$ prime to $p$. We thus have
    \[
    \Q[G]\simeq \Q[H][m,f]
    \]
    for some $f\in \Aut \Q[H]$ with $f^m=1$. By Proposition~1.4 of \cite{Yam74} the center of a simple factor of $\Q[H]$ is given by $\Q(\chi)$ for some character $\chi$ of $H$. Now, since $H$ is a $p$-group, this gives that the centers of the simple factors of $\Q[H]$ are subfields of the CM field $\Q(\mu_{p^n })$ for some $n\in \N$. As $\Q[H]$ is equipped with a positive involution, stable by $f$ -- see paragraph~\ref{subsub:twisted1} -- we are in position to apply Theorem~\ref{theo:centresalgebres}. If $p=2$ then $m$ is odd and we are in case (1). If $p$ is odd, by Satz~2 (and the following paragraph) of \cite{Roq58}, the simple factors of $\Q[H]$ are matrix algebras over their center and by Theorem~1.(i) of \cite{Fo87} these centers are of the form $\Q(\mu_{p^e})$ for some $e\geq 0$. More precisely, except for the factor $\Q$ coming from the trivial representation we have $e\geq 1$. Indeed, let $\mathrm{M}_n(\Q(\chi))=\mathrm{M}_n(\Q)$ be one of the factor with $n\geq 1$, for an irreducible non trivial character $\chi$ of $H$. By Ford, there is a linear character $\lambda$ of a subgroup of $H$ that induces $\chi$. Since $\Q(\chi)=\Q(\lambda)$ we have that $\lambda$ is trivial. But then $\chi$ can not be irreducible and non trivial as a character induced by the trivial character of a proper subgroup is a permutation character, and contains the trivial character. We are thus either in the case (2) of Theorem~\ref{theo:centresalgebres} or considering the trivial action on the factor $\Q$ corresponding to the trivial representation which only yields copies of $\Q$.

    

    \medskip

    To prove the last part of the statement we argue as the proof of Théorème~3 in \cite{Se60}. 

    \smallskip

    Let $G'\subset G$ be a subgroup satisfying $(B_q)$ for some prime number $q$. If $q=p$ then $G'$ is a $p$-group with $p$ odd so that the results of Ford and Roquette ensure that no algebra of type III appears as a simple factor of $\Q[G']$. If $q\neq p$, since $\Card G$ is odd we must have $q\neq 2$ and Lemma~\ref{lem:notypeIIIoddorder} ensures that there are not type III algebra in the decomposition of $\Q[G']$. 

    \smallskip

    Let us remark now that for a $\Q$-algebra $E$ of finite dimension with positive involution we have $E\otimes_{\Q} \R$ is quasi-split if and only if $E$ has no factor of type III. In particular, we see that for all subgroups satisfying $(B_q)$, for some prime $q$, of $G$, the algebra $\Q[G']$ is quasi-split over $\R$. By Brauer's induction theorem -- see for example Proposition~2 of \cite{Se60} -- the algebra $\Q[G]$ is also quasi-split over $\R$ which is equivalent to the fact that there is no factor of type III in its decomposition. 
\end{proof}


\section{The $(p,t,a)$-inertial groups as finite monodromy groups} \label{sec:ptaintomonod}

In this section $p$ is a fixed prime number. The goal is to apply the results of the two previous part to $(p,t,a)$-inertial groups and show that these groups are realized as finite monodromy groups of abelian varieties of dimension $g=t+a$ over number fields, answering a question of Silverberg and Zarhin. 

\subsection{Symplectic realization of finite groups and a characterization of $(p,t,a)$-inertial groups}

\subsubsection{} We start by a lemma on rational representations of simple polarizable algebras over the fields $\Q_{\ell}$. For any field extension $K/\Q$ and a $\Q$-algebra $E$ we extend the definition of a representation, a rational representation and a polarized representation of $E$ to be a $K$-vector space with an action of $E$ which is $\Q$-linear and having the corresponding properties. We are only concerned with the fields $\Q_{\ell}$ for $\ell\neq p$. To that end, let us fix an embedding $\Q_{\ell}\subset \C$ for every $\ell\neq p$. For a simple algebra $E$ quasi-split outside $p$, the results of Section~\ref{sub:ratetpolrep} extend without difficulty. 

\begin{lem} \label{lem:Qlrepratetpol}
    Let $E$ be a simple $\Q$-algebra quasi-split outside $p$ and $\ell\neq p$ a prime number. If $V$ is a rational $\C$-representation of $E$ then it is of the form $V_{\ell} \otimes_{\Q_{\ell}} \C$ for a unique $\Q_{\ell}$-representation $V_{\ell}$ of $E$ which is rational. In particular, the rational $\Q_{\ell}$-representation of $E$ are the multiples of $V_r(E)_{\ell}$ and the rational and polarized $\Q_{\ell}$-representations of $E$ are the multiples of $V_{rp}(E)_{\ell}$. 
\end{lem}

\begin{proof}
    By Theorem~\ref{theo:repratetpol} it is enough to show that the $\C$-representation $V_r(E)$ verifies the statement. 

    \medskip

    Let us write $E\simeq \mathrm{M}_n(D)$ with $D$ a skew-field of center $F$ over $\Q$ and $\deg E=ndf$ by setting $d^2=[D:F]$ and $f=[F:\Q]$. We have a decomposition
    \[
    F\otimes_{\Q} \Q_{\ell} =\prod\limits_{i=1}^r F_i
    \]
    and since $E$ is quasi-split outside $p$ we have
    \[
    E\otimes_{\Q} \Q_{\ell} \simeq \prod\limits_{i=1}^r \mathrm{M}_{nd}(F_i).
    \]
    For $i\in \{1,\dots, r\}$ let $E_i=F_i^{nd}$ be the unique non trivial simple $\mathrm{M}_{nd}(F_i)$-module. We consider the rational $\Q_{\ell}$-representation 
    \[
    V_r(E)_{\ell}= \bigoplus\limits_{i=1}^r E_i.
    \]
    We have
    \[
    V_r(E)_{\ell} \otimes_{\Q_{\ell}} \C=\bigoplus\limits_{i=1}^r (\C^{nd})^{[F_i:\Q_{\ell}]}=\bigoplus\limits_{i=1}^f \C^{nd}.
    \]
    Hence $V_r(E)_{\ell}\otimes_{\Q_{\ell}} \C$ is a rational $\C$-representation of $E$ of dimension $\deg E$. It is thus isomorphic to $V_r(E)$ which concludes.

    \medskip

    The last part of the statement is clear for rational $\Q_{\ell}$-representations. For a rational and polarized $\Q_{\ell}$-representation $V$ it is also clear that $V\otimes_{\Q_{\ell}} \C$ is rational and polarized so that $V=V_{rp}(E)_{\ell}^n$ for some integer $n\in \N$. It thus suffices to see that $V_{rp}(E)_{\ell}$ is polarized. If $E$ is of type I, II and IV this follows exactly as in the case of $\C$-representations. For $E$ of type III we thus need to show that $V_r(E)_{\ell}$ is polarized. From the proof of the first part we know that
    \[
    V_r(E)_{\ell}= \bigoplus\limits_{i=1}^r  E_i
    \]
    and the sum runs over all simple $\Q_{\ell}$-representations of $E$. In particular, we have $V_r(E)_{\ell}\simeq V_r(E)_{\ell}^{\vee}$. This isomorphism induces an isomorphism $V_r(E)\simeq V_r(E)^{\vee}$ and the induced non degenerate bilinear form on $V_r(E)$ is invariant. Since there is only one and it is symplectic we are done.  
\end{proof}

\subsubsection{} We will give a characterization of $(p,t,a)$-inertial groups through representations of algebras but first we introduce and characterize the notion of a symplectic realization of a ramification group. This will add clarity to the fact that we will not only realize the groups themselves but also any choice of a representation that satisfies part $(ii)$ of Definition~\ref{def:pta}.

\begin{defin}
    Let $G$ be a ramification group. A symplectic realization of $G$ of dimension $2a$ is the data of a quotient $E$ of $\Q[G]$ and a rational and polarized $\C$-representation $V$ of $E$ of dimension $2a$.   
\end{defin}

Note that in the definition, we consider $E$ equipped with the positive involution induced by the natural one on $\Q[G]$ and with a map $G\to E$.

\begin{theo} \label{theo:realsymp}
    Let $G$ be a ramification group. The symplectic realizations $(E,V)$ of $G$ of dimension $2a$ are in correspondence with family of maps $(p_{\ell}\colon G\to \Sp_{2a}(\Q_{\ell}))_{\ell\neq p}$ such that the characteristic polynomials of the elements of $G$ have integer coefficients independent of $\ell$. The correspondence is such that, for $\ell\neq p$, the map $p_{\ell}\colon G\to \Sp_{2a} (\Q_{\ell})$ is obtained by composition by the canonical map $G\to E$, that is $p_{\ell}$ coincides with
    \[
    G\to E \to \End_{\Q_{\ell}} V_{\ell}. 
    \]
\end{theo}
\begin{proof}
    Let $(p_{\ell}\colon G\to \Sp_{2a}(\Q_{\ell}))_{\ell\neq p}$ be a family of maps as in the statement. Let us fix $\ell\neq p$ a prime and consider $E$ the $\Q$-algebra generated by $G$ in $M_{2a}(\Q_{\ell})$. The algebra $E$ is a quotient of $\Q[G]$ by construction and one can note that it does not depend on the choice of $\ell$. It is also by construction that $E$ acts on $V_{\ell}=\Q_{\ell}^{2a}$ and this action has an invariant symplectic form. In other words, $V_{\ell}$ is a rational and polarized $\Q_{\ell}$-representation of $E$. By Lemma~\ref{lem:Qlrepratetpol} the space $V_{\ell}\otimes_{\Q_{\ell}} \C$ is a $\C$-rational and polarized representation of $E$. 

    \medskip

    Let $E$ be a quotient of $\Q[G]$ and $V$ a rational and polarized $\C$-representation of $E$. Let $\ell \neq p$ be a prime. Again, by Lemma~\ref{lem:Qlrepratetpol}, to the representation $V$ corresponds a rational and polarized $\Q_{\ell}$-representation $V_{\ell}$. Now, since $E$ is a quotient of $\Q[G]$ it is equipped with a map $\rho\colon G\to E$ and the positive involution on $E$ is such that $\rho(g)^*=\rho(g^{-1})$. In particular, since $V_{\ell}$ is polarized we see that $\rho(g)$ is an element of the symplectic group. In other words, the restriction to the image of $G$ in $E$ of the map $E\to \End_{\Q_{\ell}} V_{\ell}$ has image in $\Sp_{2a}(\Q_{\ell})$. We recover this way a family of maps $(p_{\ell}\colon G\to \Sp_{2a}(\Q_{\ell}))_{\ell\neq p}$ which satisfy the condition.

    \medskip

    It is clear that both constructions are inverse to each other. 
\end{proof}

\begin{cor}
    Let $G$ be a ramification group at $p$ and $t,a\in \N$. Then $G$ is $(p,t,a)$-inertial if and only if there is a symplectic realization $(E,V)$ of $G$ and a map $G\to \GL_t(\Z)$ such that the induced map
    \[
    G\to \GL_t(\Z)\times E
    \]
    is injective.
\end{cor}

\begin{proof}
    By definition, if $G$ is $(p,t,a)$-inertial there is a family of injective maps $(\iota_{\ell}\colon G\to \GL_t(\Z)\times \Sp_{2a}(\Q_{\ell}))_{\ell \neq p}$ for which the projections $p_{\ell}$ on the second factor satisfy the condition of Theorem~\ref{theo:realsymp}. So there is a quotient $E$ of $\Q[G]$ with a polarized $V$ representation of $E$ such that $p_{\ell}$ is obtained by the composition
    \[
    G\to E \to \End_{\Q_{\ell}}V_{\ell}.
    \]
    Since the second map is injective we have that the induced map
    \[
    G\to \GL_t(\Z) \times E
    \]
    is also injective.

    \medskip

    The converse is straightforward from Theorem~\ref{theo:realsymp}. The data of $(E,V)$ gives a family of maps $(p_{\ell}\colon G\to \Sp_{2a}(\Q_{\ell}))_{\ell\neq p}$ which gives an injection
    \[
    \iota_{\ell} \colon G\to \GL_t(\Z) \times \Sp_{2a} (\Q_{\ell})
    \]
    for all $\ell\neq p$ by construction. 
\end{proof}

\subsubsection{} We finish this section by relating symplectic realizations to geometric embeddings introduced in Section~\ref{sub:geomalgembed}. We show that all realizations are geometric in the sense of the following definition, which is the first form of our main result and the main step towards realizing the $(p,t,a)$-inertial groups as finite monodromy groups. To that end, we recall some results of Tate on abelian varieties and their Tate modules. 

\medskip

Let $A$ be an abelian variety of dimension $a$ over $\overline{\F_p}$. For a prime $\ell\neq p$, the endomorphism ring $\End A$ acts on the $\ell$-adic Tate module $T_{\ell} A=\varprojlim\limits_{n} A[\ell^n]$. Let $V_{\ell} (A)$ denote the $\Q_{\ell}$-vector space $T_{\ell} A\otimes_{\Z_{\ell}} \Q_{\ell}$. The representation of the $\Q$-algebra $E(A)=\End A \otimes_{\Z} \Q$ given by $V_{\ell}(A)$ is rational and polarized -- see \cite{Mu08} Theorem~4 on p. 167 (181). By Lemma~\ref{lem:Qlrepratetpol} the representation $V_{\ell} (A)$ corresponds to a $\C$-representation $V(A)$ of $E(A)$ and since the characteristic polynomials of the representation $V_{\ell} (A)$ are independent of $\ell$ the $\C$-representation $V(A)$ is independent of $\ell$. For an abelian variety $A$ over a finite field, we will keep the notation $E(A)$ and $V(A)$ for its endomorphism algebra and the $\C$-representation of $E(A)$ associated to its $\ell$-adic Tate module. Note that when $E(A)$ is simple, by the results of Tate on the dimension of $E(A)$ we have in this case that $V(A)=V_{rp}(E(A))$. 

\begin{defin}
    A symplectic realization $(E,V)$ of a ramification group $G$ is said to be geometric if there is an abelian variety $A$ over $\overline{\F_p}$ and an embedding $E\subset E(A)$ such that the induced $\C$-representation of $E$ given by $V(A)$ is isomorphic to $V$. 
\end{defin}

\begin{theo} \label{theo:realgeomall}
    Let $G$ be a ramification group at $p$. All symplectic realizations of $G$ are geometric. 
\end{theo}

\begin{proof}
Let $(E,V)$ be a symplectic realization of $G$. First note that if $E$ is not simple and we have a decomposition in simple factors $E\simeq \prod_{i=1}^r E_i$ then $V$ splits in a direct sum of representations $V\simeq \bigoplus_{i=1}^r V_{rp}(E_i)^{n_i}$ for some $n_i\in \N$. For each $i\in \{1,\dots, r\}$ we get a symplectic representation $(E_i, V_{rp}(E_i)^{n_i})$ of $G$ and it is clear that $(E,V)$ is geometric if and only if $(E_i,  V_{rp}(E_i)^{n_i})$ is geometric for each $i\in \{1,\dots, r\}$. We thus assume $E$ is simple. 

\medskip

We start by showing that $(E,V_{rp}(E))$ is a symplectic realization of $G$ which is geometric. Since $G$ is a ramification group, by $(i)$ of Theorem~\ref{theo:structureQG} the algebra $E$ is quasi-split outside $p$. If $E$ is of type I, II or IV then by Theorem~\ref{theo:goodembed} we have a good embedding $E\to E(A)$ for some abelian variety $A$ over $\overline{\F_p}$. We thus have that $V(A)=V_{rp}(E(A))=V_{rp}(E)$ and the realization $(E,V)$ is geometric. 

\smallskip

Let us now consider the case where $E$ is of type III. By $(ii)$ of Theorem~\ref{theo:structureQG} the center of $E$ is subfield of $\Q(\mu_{p^{n}})$ for some $n\in \N$. But this field is totally ramified at $p$ and thus also the center $F$ of $E$. Hence there is only one Brauer class of type III algebras over $F$, the algebras of the form $\mathrm{M}_k(F\otimes \HH_{p,\infty})$. It follows that $E$ is of that form and by Theorem~\ref{theo:goodembed} again there is $A$ over $\overline{\F_p}$ with $V(A)=V_{rp}(A)=V_{rp}(E)$.   

\medskip

We go back to the general case, that is of a realization $(E,V)$ of $G$ with $E$ simple. By Theorem~\ref{theo:repratetpol} the representation $V$ is isomorphic to $V_{rp}(E)^n$ for some $n\in \N$. From the previous case, let $B$ be an abelian variety over $\overline{\F_p}$ with $E\hookrightarrow E(B)$ such that $V(B)=V_{rp}(E)$. The abelian variety $A=B^n$ suits our needs. Indeed, we have a diagonal embedding $E\hookrightarrow E(A)\simeq \mathrm{M}_n(E(B))$ for which the induced $\C$-representation $V(A)$ has the desired property $V(A)=V_{rp}(E)^n=V$ as $\C$-representations of $E$. 
\end{proof}

\subsection{The realization of $(p,t,a)$-inertial groups as finite monodromy groups}

\subsubsection{} Let us start by recalling the results of \cite{Ph4} which relate to finite monodromy groups. For an abelian variety $A$ over a $p$-adic field $K$, the finite monodromy group $\Phi(A)$ of $A$ is the Galois group of the smallest field extension $L_A$ of the maximal unramified extension $K^{\mathrm{un}}$ of $K$ such that $A_{L_A}$ has semi-stable reduction -- the fact that $L_A/K^{\mathrm{un}}$ is Galois is section 4 of Exposé IX of \cite{sga}. 

\medskip

In \cite{Ph4}, based on chapter 2 of \cite{FC90}, there is a notion of a polarized semi-abelian variety. We recall it here as it is used in the main theorem. 

\begin{defin} Let $A_0$ and $A_0^t$ be semi-abelian varieties over a finite field. A polarization $\lambda_0\colon A_0\to A_0^t$ is a map of semi-abelian varieties
\[
\begin{tikzcd}
    0 \arrow[r] & T_0 \arrow[r] \arrow[d, "\lambda_{T_0}"] & A_0 \arrow[r] \arrow[d, "\lambda_0"] & B_0 \arrow[r]  \arrow[d, "\lambda_{B_0}"]& 0\\
    0 \arrow[r] & T_0^t \arrow[r] & A_0^t \arrow[r] & B_0^{\vee} \arrow[r] & 0
\end{tikzcd}
\]
such that the induced maps $\lambda_{T_0}$ is an isogeny and $\lambda_{B_0}$ is a polarization. A semi-abelian variety $A_0$ equipped with a polarization $\lambda_0$ is said to be polarized.
    
\end{defin}

In the next sections of the text, we will conserve the notations $T_0$ and $B_0$ for the maximal torus and abelian quotient of a semi-abelian variety $A_0$. Note that, for a prime $\ell\neq p$, through the $\ell$-adic Tate module of $B_0$ there is a map 
\[
\Aut A_0\longrightarrow \GL_t(\Z)\times \Sp_{2a} (\Q_{\ell})
\]
where $t$ is the dimension of $T_0$ and $a$ the dimension of $B_0$.

\medskip

The main theorem of \cite{Ph4} now states that a ramification group $G$ is the finite monodromy group of some abelian variety if and only if there is an embedding $G\subset \Aut (A_0, \lambda_0)$ for some polarized semi-abelian variety $A_0$ over $\overline{\F_p}$ -- see Théorème~1.1 of \cite{Ph4}. 

\subsubsection{} We can now show that all $(p,t,a)$-inertial groups are realized as finite monodromy groups in dimension $t+a$. More precisely, we show that for a ramification group $G$ all symplectic realizations are obtained by realizations of $G$ as a finite monodromy group. We first give the context for a precise meaning of this last part. For an abelian variety $A$ over a $p$-adic field $K$, its finite monodromy group $\Phi_A$ acts on the reduction of $A_{L_A}$, which is a semi-abelian variety $A_0$ over $\overline{\F_p}$. This action produces an injection $G\hookrightarrow \Aut (A_0,\lambda_0)$ and a symplectic realization $(E(B_0),V(B_0))$. We show that all symplectic realizations of $G$ are obtained in this way. 

\begin{theo} \label{theo:realmonod}
    Let $G$ be a ramification group at $p$. Then for all family of injective maps $(\iota_{\ell}\colon G\to \GL_t(\Z)\times \Sp_{2a}(\Q_{\ell}))_{\ell\neq p}$ with second projection corresponding to a symplectic realization $(E,V)$ of dimension $2a$ of $G$ there is an abelian variety $A$ of dimension $t+a$ over a $p$-adic field $K$ such that $G$ is the finite monodromy group of $A$. Moreover, the action of $G$ on the reduction $A_0$ of $A_{L_A}$ recover the first projection $G\to \GL_t(\Z)$ of $\iota_{\ell}$ for any $\ell\neq p$ and the symplectic realization $(E,V)$ up to isomorphism. 
\end{theo}

\begin{proof}
    Let $(E,V)$ be the symplectic realization coming from the second projection of the family of maps $(\iota_{\ell})_{\ell\neq p}$. From Theorem~\ref{theo:realgeomall} this realization is geometric so there is an abelian variety $B_0$ over $\overline{\F_p}$ such that $E\subset E(B'_0)$ and $V(B'_0)=V$. Consider the image of the order $\Z[G]$ of $\Q[G]$ in $E(B'_0)$. It is contained in a maximal order and by Theorem~3.13 of \cite{Wa69} there is an abelian variety $B_0$ for which $E(B_0)=E(B'_0)$ and $V(B_0)=V(B'_0)$ but also such that $\End B_0\subset E(B_0)$ contains the image of $G$. Since $G$ is finite there is a polarization $\lambda_0$ of $B_0$ which is stable by $G$ and we can consider the semi-abelian variety $A_0=\G_m^t\times B_0$ over $\overline{\F_p}$. The semi-abelian variety $A_0$ is polarized by $(\mathrm{id}, \lambda_0)$ and we have an injective map
    \[
    G\longrightarrow \Aut A_0= \GL_t(\Z) \times \Aut B_0.
    \]
    By Théorème~1.1 of \cite{Ph4} we get an abelian variety $A$ over $p$-adic field $K$ such that $G$ is the finite monodromy group of $A$ and the reduction of $A_{L_A}$ is isogenous to $A_0$. By construction, the symplectic realization obtained by $(A,G)$ is isomorphic to $(E,V)$. 
\end{proof}

\begin{cor}\label{cor:realmonod}
    Let $G$ be a ramification group. Then $G$ is $(p,t,a)$-inertial if and only if $G$ is the finite monodromy group of an abelian variety $A$ of dimension $t+a$ over a number field $K$ at a place $v$ of residue characteristic $p$. Furthermore, we can choose $A$ such that, for $L/K$ an extension for which $A_L$ has semi-stable reduction at the places of $L$ above $v$, the reduction of $A_L$ at the places above $v$ has toric rank $t$ and abelian rank $a$. 
\end{cor}

\begin{proof}
    By Theorem~\ref{theo:realmonod} we get that if $G$ is $(p,t,a)$ inertial such an abelian variety $A'$ over a $p$-adic field exists. By Theorem~4.3 of \cite{Ph221} we get $A$ over a number field which satisfies all our conditions. The converse is given, for example, by Theorem~5.2 of \cite{SZ98}.
\end{proof}

\subsubsection{} We finish by giving some examples of computations on specific groups. For a ramification group $G$, we are able to determine the values for which it is $(p,t,a)$-inertial in a methodic way with finitely many steps. Note that, from the regular representation of $G$, it is clear that there are minimal integers $t_G,a_G\geq 1$ such that $G$ is $(p,t_G,0)$ and $(p,0,a_G)$-inertial. It is also clear that if $G$ is $(p,t,a)$-inertial for some $t,a\geq 0$ then it is trivially $(p,t+n,a+m)$-inertial for any $n,m\geq 0$. We thus see that there is only finitely many couples $(t,a)$ for which $G$ is not $(p,t,a)$-inertial, inside a box of height $t_G-1$ and length $a_G-1$ starting from $(0,0)$ in the lattice $\Z^2$. Lastly, we remark that if $G$ is $(p,t,a)$-inertial with $t,a\geq 1$ in a non-trivial way, that is for fixed $\ell \neq p$ the map
\[
\iota_{\ell} \colon G\longrightarrow \GL_t(\Z) \times \Sp_{2a}(\Q_{\ell})
\]
has its images in each factor not $G$, then $G$ is a subdirect factor of the product of two of its non trivial normal subgroups. In other words, there are normal proper subgroups $K$ and $H$ in $G$ such that $K\cap H =\{1\}$ such that the natural map
\[
G\longrightarrow G/K\times G/H
\]
is injective. For a group for which no such subgroups exists, the set of exceptional couples $(t,a)$ for which it is not $(p,t,a)$-inertial will exactly be the box mentioned previously. 

\begin{ex}
    We give three examples.
    \begin{itemize}
        \item[$(i)$] Let $n\geq 1$ and $G=\Z/n\Z$. Cyclic groups are ramification groups at all primes $p$. We compute the list of exceptional couples $(t,a)$ for $G$, which does not depend on $p$. Let $D_n$ be the set of ordered couples $(r,s)$ of factors of $n$ such that $r$ is prime to $s$ and $n=r\cdot s$. Considering the factors $\Q(\mu_r)$ and $\Q(\mu_s)$ of $\Q[G]$ we get a representation $G\to \GL_{\varphi(r)}(\Z)$ with image $\Z/r\Z$ on one side and a symplectic realization $(\Q(\mu_s), V_{rp}(\Q(\mu_s)))$ on the other. This symplectic realization gives a family of maps, for all primes $\ell$, 
        \[
        p_{\ell} \colon G\to \Sp_{\varphi(s)} (\Q_{\ell})
        \]
        with images $\Z/s\Z$. The product maps 
        \[
        G\longrightarrow \GL_{\varphi(r)} (\Z) \times \Sp_{\varphi(s)}(\Q_{\ell})
        \]
        are injective by construction and make $G$ into a $(p,\varphi(r),\frac{\varphi(s)}{2})$-inertial group. It is an easy check that for any value $t,a\geq 0$ with $t\leq \varphi(r)-1$ or $a\leq \frac{\varphi(s)}{2}-1$ for any couple $(r,s)\in D_n$, $G$ is not $(p,t,a)$-inertial. 

        \smallskip

        We thus recover and generalize the result of Silverberg and Zarhin for small cyclic groups in \cite{SZ05}. 

        \item[$(ii)$] We consider the wreath product $G=Q_8\wr \Z/2\Z$. This group is shown to be $(2,0,2)$-inertial by Silverberg and Zarhin in \cite{SZ05} and was realized as a finite monodromy group over number fields by Chrétien and Matignon in \cite{CM13}. Note that in this group, all pair of normal subgroups intersect non-trivially so that it can not be a $(p,t,a)$-inertial group with $t,a\geq 1$ in a non-trivial way. We thus look for the smallest rational and polarized faithful representation coming from a factor of its rational algebra in order to compute $t_G$ and $a_G$. Let us start by computing the rational group algebra of $Q_8\times Q_8$. We have
        \begin{align*}
        \Q[Q_8\times Q_8] &= ( \Q\times \Q \times \Q \times \Q\times \HH_{2,\infty}) \otimes_{\Q} ( \Q\times \Q \times \Q \times \Q \times \HH_{2,\infty}) \\
          &=  \Q^4\times (\Q^2)^{6} \times (\HH_{2,\infty}^2)^4 \times \mathrm{M}_{4}(\Q)
        \end{align*}
        where we regrouped the factors that are swapped by the $\Z/2\Z$-action coming from $G$. The factors that are fixed get doubled in $\Q[G]$ and those that are swapped give $2\times 2$ matrix algebra over the respective division algebra by Lemma~\ref{lem:isomcentrealgebreint}. We thus have
        \[
        \Q[G]= \Q^8\times \mathrm{M}_{2}(\Q)^6 \times \mathrm{M}_4(\Q)^2 \times \mathrm{M}_2(\HH_{2,\infty})^4.
        \]
        Now, it is an easy check that we have an embedding $G\hookrightarrow \mathrm{M}_2(\HH_{2,\infty})$ so that $G$ is both $(2,0,2)$-inertial since $\deg \mathrm{M}_2(\HH_{2,\infty})=4$ and it is a type III algebra, and $(2,8,0)$-inertial since the Schur index of $\HH_{2,\infty}$ is $2$. We thus have $t_G=8$ and $a_G=2$ and since $G$ is not a subdirect product it is not $(2,t,a)$-inertial for any $(t,a)$ in the box delimited by $(0,0)$, $(7,0)$, $(7,1)$ and $(0,1)$ except for the values in the top and right edges.

        \item[$(iii)$] Consider the dihedral group $G=D_{13}$ of order $26$. By the table on page 105 of \cite{VM02} we have
        \[
        \Q[G]=\Q\times \Q \times \mathrm{M}_2(\Q(\zeta_{13}+\zeta_{13}^{-1}).
        \]
        We get a faithful representation from the factor $\mathrm{M}_2(\Q(\zeta_{13}+\zeta_{13}^{-1}))$ which yield by a simple computation $t_G=6$ and $a_G=3$. Since $G$ is not a subdirect product again, the exceptional values are in the box delimited by $(0,0)$, $(5,0)$, $(5,2)$ and $(0,2)$.
    \end{itemize}
\end{ex}

In general, given a ramification group $G$ one can determine the values for which it is $(p,t,a)$-inertial by computing the rational group algebra $\Q[G]$ and the faithful representations given by the quotients of $\Q[G]$. Then, the values are computed using the dimension of the given representations and the degrees of the polarized representation of the simple algebras in play.

\printbibliography

\end{document}